\newif\iffinal
\finaltrue	% Final version

\documentclass[letterpaper,11pt,reqno]{amsart} 
\RequirePackage[OT1]{fontenc}
\usepackage[portrait,margin=3.5cm]{geometry}  
\iffinal\else\usepackage[notref,notcite]{showkeys}\fi
\usepackage{mathrsfs}
\usepackage[colorlinks,citecolor=blue,urlcolor=blue,linkcolor=blue,pagecolor=blue,linktocpage=true,backref=true]{hyperref}
\usepackage[foot]{amsaddr}
\usepackage{amssymb,amsthm,amsfonts,amsbsy,latexsym,dsfont}
\usepackage[textsize=small]{todonotes}       % includes TO DO LIST
\usepackage{graphicx}
\usepackage[numeric,initials,nobysame]{amsrefs}
\usepackage{upref,setspace,bbm}
\usepackage{todonotes}

\usepackage{enumerate}

\newenvironment{enumeratea}{\begin{enumerate}[\upshape (a)]}{\end{enumerate}}

\numberwithin{equation}{section}
\numberwithin{figure}{section}
\numberwithin{table}{section}

\newtheorem{thm}{Theorem}[section]
\newtheorem{theorem}{Theorem}[section]
\newtheorem{lem}[thm]{Lemma}
\newtheorem{cor}[thm]{Corollary}
\newtheorem{prop}[thm]{Proposition}

\renewcommand{\leq}{\leqslant} 
\renewcommand{\geq}{\geqslant}

\newcommand{\ind}{\mathds{1}}
\newcommand{\eps}{\varepsilon}

\newcommand{\set}[1]{\left\{#1\right\}}

\newcommand{\Probc}{\stackrel{\mathrm{P}}{\longrightarrow}}
\newcommand{\weakc}{\stackrel{\mathrm{w}}{\longrightarrow}}

\def\qed{ \hfill $\blacksquare$}  
\newcommand{\sss}{\scriptscriptstyle}
\newcommand{\op}{o_{\sss \pr}}
\newcommand{\Op}{O_{\sss \pr}}
\newcommand{\eqn}[1]{\begin{equation}#1\end{equation}}
\newcommand{\eqan}[1]{\begin{align}#1\end{align}}
\newcommand{\e}{{\mathrm e}}
\newcommand{\nn}{\nonumber}

\def\ind{{\rm 1\hspace{-0.90ex}1}}
\newcommand\1{\mathbbm{1}}
\newcommand{\indic}[1]{\1_{\{#1\}}}

\newcommand{\Flood}{{\mathrm{Flood}}}
\newcommand{\diam}{{\mathrm{Diam}}}

\newcommand{\SWG}{{\sf SWG}}
\newcommand{\cA}{\mathcal{A}}\newcommand{\cB}{\mathcal{B}}
\newcommand{\cD}{\mathcal{D}}\newcommand{\cE}{\mathcal{E}}\newcommand{\cF}{\mathcal{F}}
\newcommand{\cG}{\mathcal{G}}\newcommand{\cH}{\mathcal{H}}
\newcommand{\cK}{\mathcal{K}}\newcommand{\cL}{\mathcal{L}}
\newcommand{\cN}{\mathcal{N}}
\newcommand{\cP}{\mathcal{P}}\newcommand{\cR}{\mathcal{R}}

\newcommand{\cV}{\mathcal{V}}
  
\newcommand{\mvD}{\boldsymbol{D}}

\newcommand{\mvd}{\boldsymbol{d}}

\newcommand{\mvpi}{\boldsymbol{\pi}}

\newcommand{\bL}{\mathbb{L}}

\newcommand{\bR}{\mathbb{R}}

\newcommand{\bZ}{\mathbb{Z}}        
\newcommand{\sS}{\mathscr{S}}

\DeclareMathOperator{\E}{\mathds{E}}
\DeclareMathOperator{\pr}{\mathds{P}}
\DeclareMathOperator{\prob}{\mathds{P}}

\DeclareMathOperator{\var}{Var}

\newcommand{\tra}{\mbox{tr}}
\newcommand{\tca}{\mbox{co}}

\definecolor{darkgreen}{rgb}{0,.4,0}
\definecolor{darkagenta}{rgb}{.5,0,.5}
\definecolor{darkred}{rgb}{1,0,0}%was 0.85
\definecolor{darkblue}{rgb}{0,0,.4}

\newcommand{\ch}[1]{#1}
\newcommand{\csb}[1]{#1}
\newcommand{\crh}[1]{#1}

\begin{document}

\title[Diameter of the stochastic mean-field model]{Diameter of the stochastic\\
mean-field model of distance}

\date{\today}

\author[Bhamidi]{Shankar Bhamidi$^1$}
\address{$^1$Department of Statistics, University of North Carolina, Chapel Hill,}
\author[van der Hofstad]{Remco van der Hofstad$^2$}
\address{$^2$Department of Mathematics and
    Computer Science, Eindhoven University of Technology, P.O.\ Box 513,
    5600 MB Eindhoven, The Netherlands.}
\email{bhamidi@email.unc.edu, rhofstad@win.tue.nl}

\maketitle
\begin{abstract}
We consider the complete graph $\cK_n$ on $n$ vertices with exponential mean $n$ edge lengths. Writing $C_{ij}$ for the \crh{weight of the smallest-weight} path between vertex $i,j\in [n]$, Janson \cite{janson-cm-gp} showed \ch{that $\max_{i,j\in [n]} C_{ij}/\log{n}$ converges in probability to 3. We extend this results by showing that 
$\max_{i,j\in [n]} C_{ij} - 3\log{n}$ converges in distribution to some limiting random variable that} can be identified via a maximization procedure on a limiting infinite random structure.  Interestingly, this limiting random variable has also appeared as the weak limit of the re-centered \emph{graph diameter}   
of the barely supercritical Erd\H{o}s-R\'enyi random graph in \cite{RioWor10}.
\end{abstract}

%%%%%%%%%%%%%%%%%%%%%%%%%%%%%%%%%%%%
\section{Introduction}
\label{sec:int}
We consider the complete graph $\cK_n$ on the vertex set $[n]:=\set{1,2,\ldots, n}$ and edge set $\cE_n :=\set{\set{i,j}: i< j \in [n]}$. To each edge $e\in \cE_n$, assign exponential mean $n$ edge lengths $E_e$, independently across edges. This implies for any vertex $v$, the closest neighbor to this vertex is $O_P(1)$ distance away. Define the length of a path $\pi$ as
	\begin{equation}
	\label{eqn:wt-path-def}
	w(\pi):= \sum_{e\in \pi} E_e.
	\end{equation} 
This assignment of random edge lengths makes $\cK_n$ a (random) metric space often referred to as the \emph{stochastic mean-field model of distance} (see Section \ref{sec:disc}). By continuity of the distribution of edge lengths, this metric space has unique geodesics. For any two vertices $i,j\in [n]$, let $\mvpi(i,j)$ denote the shortest path between these two vertices and write $C_{ij}$ for the length of this geodesic. The functional of interest in this paper is the diameter of the metric space:
	\begin{equation}
	\label{eqn:diam-defn}
	\diam_w(\cK_n):=\max_{i,j\in [n]} C_{ij}.
	\end{equation}
We first dive into the statement of the main result, postponing a full discussion to Section \ref{sec:disc}.

\section{Results}
\label{sec:res}
The main aim of this paper is to prove that the diameter defined in \eqref{eqn:diam-defn} properly re-centered converges to a limiting random variable. We start by constructing this \ch{limiting} random variable.

\subsection*{Construction of the limiting \ch{random variable}}
The \ch{limiting} random variable arises as an optimization problem on an infinite randomly weighted graph $\cG_\infty = (\cV,\cE)$.  The vertex set of this graph is the set of positive integers $\bZ_+ = \set{1,2,\ldots}$, while the edge set consists of all undirected edges $\cE = \set{\set{i,j}: i,j\in \bZ_+, i\neq j}$. Let $\cP$ be a Poisson process on $\bR$ with intensity measure having density 
	\begin{equation}
	\label{eqn:pp-rate-def}
	\lambda(y) = \e^{-y}, \qquad -\infty < y < \infty.
	\end{equation}
It is easy to check that $\max\set{x: x\in \cP } < \infty$ a.s. Thus we can order the points in $\cP$ as $Y_1> Y_2> \cdots$. We think of $Y_i$ as the vertex weight at $i\in \bZ_+$. The edge weights are easier to describe.  Let $(\Lambda_{st})_{s,t\in \bZ_+, s< t}$ be a family of independent standard Gumbel random variables, namely $\Lambda_{st}$ has cumulative distribution function 
	\begin{equation}
	\label{eqn:gumbel-cdf}
	F(x) = \e^{-\e^{-x}}, \qquad -\infty < x < \infty. 
	\end{equation} 
The random variable $\Lambda_{s,t}$ gives the weight of an edge $\set{s,t} \in \cE$.  Now consider the optimization problem 
	\begin{equation}
	\label{eqn:xi-def}
	\Xi:= \crh{\max_{s,t\in \bZ_+, s< t} (Y_s+ Y_t - \Lambda_{st})}. 
	\end{equation}
Though not obvious, we shall show that $\Xi< \infty$ a.s. The main result in this paper is as follows. We write $\weakc$ to denote convergence in distribution.

\begin{theorem}[Diameter asymptotics]\label{thm:main}
For the diameter of the stochastic \crh{mean-field model} of distance, as $n\to\infty$ 
	\[
	\max_{i,j\in [n]} C_{ij} - 3\log{n} \weakc \Xi,
	\]
and
	\begin{equation}
	\label{eqn:mom-convg}
	\E[\max_{i,j\in [n]} C_{ij}] - 3\log{n} \to \E[\Xi],
	\qquad
	{\rm Var}(\max_{i,j\in [n]} C_{ij})\to {\rm Var}(\Xi).
	\end{equation}
\end{theorem}
\noindent{\bf Remark:} Theorem \ref{thm:main} solves \cite[Problems 1 and 2]{janson-cm-gp}.

\subsection{Basic notation}
\label{sec:not}
Let us briefly describe the notation used in the rest of the paper. We write
$\Probc$ to denote convergence in probability. For a sequence of random variables
$(X_n)_{n\geq 1}$, we write $X_n=\Op(b_n)$ when $|X_n|/b_n$
is a tight sequence of random variables as $n\rightarrow\infty$,
and $X_n=\op(b_n)$ when $|X_n|/b_n\Probc 0$ as $n\rightarrow\infty$.
For a non-negative function $n\mapsto g(n)$,
we write $f(n)=O(g(n))$ when $|f(n)|/g(n)$ is uniformly bounded, and
$f(n)=o(g(n))$ when $\lim_{n\rightarrow \infty} f(n)/g(n)=0$.
Furthermore, we write $f(n)=\Theta(g(n))$ if $f(n)=O(g(n))$ and $g(n)=O(f(n))$.
Finally, we write that a sequence of events $(A_n)_{n\geq 1}$
occurs \emph{with high probability} (whp) when $\prob(A_n)\rightarrow 1$. We 
use $Y\sim \exp(\lambda)$ to denote a random variable which has an exponential 
rate $\lambda$ distribution. 

\section{Background and related results}
\label{sec:disc}
We now discuss our results and place them in the context of results in the literature.

\subsection{Stochastic mean-field model of distance}
The stochastic mean-field model of distance has arisen in a number of different contexts in understanding the structure of combinatorial optimization problems in the presence of random data, ranging from shortest path problems \cite{janson-cm-gp}, random assignment problems \cite{aldous-asgn-1,aldous-asgn-2}, minimal spanning trees \cite{frieze1985value,janson-mst} and traveling salesman problems \cite{wastlund2010mean}; see \cite{aldous-steele} for a comprehensive survey and related literature. 
The closest work to this study is the paper by Janson \cite{janson-cm-gp}. Recall that $C_{ij}$ denotes the length of the geodesic between two vertices $i,j\in [n]$; by symmetry this has the same distribution for any two vertices in $i,j$. For any vertex $i\in [n]$, write $\Flood[i] := \max_{j\in [n]} C_{ij}$ for the maximum time started at $i$ to reach all vertices in $\cK_n$ (often called the \emph{flooding} time). Then Janson proved that as $n\to\infty$,
	\begin{equation}
	\label{eqn:first-order-janson}
	\frac{C_{ij}}{\log{n}} \Probc 1, \qquad \frac{\Flood[i]}{\log{n}} \Probc 2, \qquad 
	\frac{\diam_w(\cK_n)}{\log{n}} \Probc 3,	
	\end{equation} 
and further 
	\begin{equation}
	\label{eqn:sec-or-cij}
	C_{ij} - \log{n} \weakc \Lambda_1 + \Lambda_2- \Lambda_{12},
	\end{equation}
while 
	\begin{equation}
	\label{eqn:sec-or-flood}
	\Flood[i] - 2\log{n} \weakc \Lambda_1 + \Lambda_2.
\end{equation}
\ch{Here} $\Lambda_1, \Lambda_2, \Lambda_{12}$ are all independent standard Gumbel random variables as in \eqref{eqn:gumbel-cdf}. Problems 1 and 2 in \cite{janson-cm-gp} then ask if one \ch{expects a similar result as in \eqref{eqn:sec-or-cij} and \eqref{eqn:sec-or-flood}} for the diameter $\diam_w(\cK_n)$ (by \eqref{eqn:first-order-janson} obviously re-centered by $3\log{n}$).

The main aim of this paper is to answer\ch{this question} in the affirmative. We discuss more results about the distribution of $\Xi$ in Section \ref{sec:lim-max-prob}.  In the context of \eqref{eqn:mom-convg}, for $C_{ij}$ and $\Flood[i]$, Janson also shows convergence of the expectation and variance with explicit limit constants. We have been unable to derive explicit values for the limit constants $\E(\Xi)$ and $\var(\Xi)$.

\subsection{Hopcount and extrema}
This paper looks at the length of optimal paths (measured in terms of the edge weights).  One could also look at the \emph{hopcount} or the number of edges $|\mvpi(i,j)|$ on the optimal path as well as the longest hopcount $\cD^\star = \max_{i,j\in[n]} |\mvpi(i,j)|$. The entire shortest path tree from a vertex $i$ has the same distribution as a random recursive tree on size $n$ vertices (see \cite{smythe-mahmoud-rr} for a survey). Janson used this in \cite{janson-cm-gp} to show that
	\[
	\frac{|\mvpi(i,j)| - \log{n}}{\sqrt{\log{n}}} \weakc Z,
	\]
\ch{where $Z$ has a standard normal distribution.} The maximal hopcount $\cH_n(i) = \max_{j\in[n]}|\mvpi(i,j)|$ from a vertex $i$ has the same distribution as the height of random recursive tree, which by \cite{devroye1987branching} or \cite{Pitt94} satisfies the 
asymptotics $\cH_n(i)/\log{n} \Probc \e$ as $n\to\infty$. 

\ch{The first order} asymptotics for the maximum hopcount $\cD^\star$ were recently proved in \cite{berry-lugosi-broutin}, showing that $\cD^\star/\log{n} \Probc \alpha^\star$ where $\alpha^\star\approx 3.5911$ is the unique solution of the equation $x\log{x} - x = 1$.

\subsection{First passage percolation on random graphs}
The last few years have seen progress in the understanding of optimal paths in the presence of edge disorder (usually assumed to have exponential distribution) in the context of various random graph models (see e.g \cite{remco-gerard-fpp,amini2011diameter,amini2012shortest} and the references therein). In particular,  Proposition \ref{prop-many-vertices-weights} \ch{below} with a sketch of proof has appeared in \cite{shankar-er,david-shankar,aldous-exch}. 

In the context of our main result, \cite{amini2011diameter} studied the weighted diameter for the random $r$-regular graphs $\cG_{n,r}$ with exponential edge weights and proved first order asymptotics. We conjecture that \ch{one can adapt the main techniques in this paper to show the second order asymptotics}
for $r\geq 3$, i.e., 
	\eqn{
	\diam_w(\cG_{n,r})-\left(\frac{1}{r-2}+\frac{2}{r}\right)\log{n}\weakc \Xi_r,
	}
for a limit random variable $\Xi_r$ \ch{that satisfies that,} as $r\to \infty$,
	\eqn{
	r\Xi_r\weakc \Xi.
	}

\subsection{Diameter of the barely supercritical Erd\H{o}s-R\'enyi random graph}
Consider the barely supercritical Erd\H{o}s-R\'enyi random graph $\cG_n (n, (1+\eps)/n)$ where $\eps = \eps_n\to 0$ but $\eps n^3\to \infty$. It turns out that the random variable $\Xi$ in Theorem \ref{thm:main} is closely related to the random variable describing second order fluctuations for the {\bf graph diameter} $\diam_g(\cG_n (n, (1+\eps)/n))$. Here we use $\diam_g(\cdot)$ for the graph diameter of a graph, namely the largest {\bf graph} distance between any two vertices in the same component. We now describe this result. Consider the minor modification of the optimization problem defining $\Xi$ in Section \ref{sec:res} where the Poisson process $\cP$ generating the vertex weights has intensity measure with density 
	\[
	\lambda (y) = \gamma \e^{-y}, \qquad -\infty < y < \infty.
	\]
\ch{As before,} the edge weights $\Lambda_{st}$ are independent standard Gumbel random variables. Let $\Xi_\gamma$ denote the random variable corresponding to the optimization problem in \eqref{eqn:xi-def}. Let $\lambda= 1+\eps$ and let $\lambda_*< 1$ be the unique value satisfying $\lambda_*\e^{-\lambda_*} = \lambda \e^{-\lambda}$. After an initial analysis in \cite{DinKimLubPer09, DinKimLubPer11}, Riordan and Wormald in \cite[Theorem 5.1]{RioWor10} showed that there exists a constant $\gamma > 0$ such that 
	\[
	\diam_g(\cG_n (n, (1+\eps)/n)) - \frac{\log{\eps^3 n}}{\log{\lambda}}  
	- 2\frac{\log{\eps^3 n}}{\log{1/\lambda_*}} \weakc \Xi_\gamma.
	\]
We believe that the Poisson cloning technique in \cite{DinKimLubPer09, DinKimLubPer11} coupled with the techniques in this paper \ch{may yield} an alternate proof of this result but we defer this to future work.

\section{Proofs}
\label{sec:proofs}
We start with the basic ideas behind the main result. We then describe the organization of the rest of the section which deals with converting this intuitive picture into \ch{proper} proof. 

\subsection{Proof idea}
\label{sec:proof-idea}
%For the rest of the proof, we rescale edge lengths by a factor $n$ and thus assume that each edge in $\cK_n$ has an exponential mean $n$ edge length. This only serves to rescale all distances by a factor $n$ and thus reduces notational overhead. 
We write $\sS_n = (\cK_n, \set{E_e: e\in \cE_n})$ for the (random) metric space where $(E_e)_{e\in \cE_n}$ are i.i.d.\ mean $n$ exponential random variables. 
Now note that by Janson's result (\eqref{eqn:sec-or-cij}), the distance $C_{ij}$ between typical vertices $i,j\in [n]$ scales like $\log{n}+O_P(1)$. Intuitively, the extra $2\log{n}$ in the diameter arises due to the following reason. Consider ranking the vertices according to the distance to  their closest neighbor. More precisely, for each vertex $i\in [n]$, write $X_{\sss(i)} = \min_{j\in [n], j\neq i} E_{ij}$, the distance to the closest vertex to $i$. Arrange these as $X_{\sss(V_1)} > X_{\sss(V_2)} > \cdots X_{\sss(V_n)}$.  We shall show that: 
 	
	\begin{enumeratea}
 	\item the point process $\cP_n = (X_{\sss (V_i)}-\log{n}: i\geq 1)$ converges to the Poisson point process $\cP$ in Section \ref{sec:res} with intensity measure given by \eqref{eqn:pp-rate-def};
    	\item the diameter of $\cK_n$ corresponds to the shortest path between a pair of these ``slow'' vertices $(V_s, V_t)$;
	\item further, after reaching the closest vertex, the remaining path behaves like a typical optimum path in the original graph $\cK_n$ equipped with exponential mean $n$ edge lengths, but now between 2 disjoint pairs of vertices.
 \end{enumeratea}  

\ch{More precisely, part (c) entails that} $C_{V_s, V_t}  \approx X_{\sss(V_s)} + X_{\sss(V_t)} + d_w(A, B)$ where $A= \set{a,b}$ \ch{with $a,b,c,d$ four distinct vertices in $[n]$} and $d_w(A,B)$ is a random variable independent of $X_{\sss(V_t)}, X_{\sss (V_s)}$ having the same distribution as the distance between the sets $A,B$ in the original metric space $\sS_n$.  The first two terms correspond to the time to get out of these ``slow'' vertices, \ch{which scale} like $\log{n}+O_P(1)$ by (a) while $d_w(A,B)$ scales like $\log{n}+ O_P(1)$, thus implying that the diameter scales like $3\log{n}+ O_P(1)$.  \ch{By investigating the fluctuations of 
$X_{\sss(V_s)}, X_{\sss(V_t)}$ and $d_w(A, B)$, we can also identify the fluctuations of $n \max_{i,j\in [n]} C_{ij}$.} 
\bigskip

\noindent{\bf Organization of the proof:} 
 We start in Section \ref{sec:explicit} by describing the distribution of the shortest path between two disjoint set of vertices. Section \ref{sec:pp-limit} proves a weaker version of the Poisson point process limit described in (a) above. Section \ref{sec:jt-con} describes the limiting joint distribution of the (properly re-centered) weights of optimal paths between multiple source destination pairs in $\sS_n:= (\cK_n, \set{E_e: e\in \cE_n})$. Section \ref{sec:dis-min-mult} uses the results in Section \ref{sec:pp-limit} and \ref{sec:jt-con} to study asymptotics for the joint distribution of distances between the slow vertices $(V_s)_{s\in [n]}$. Section \ref{sec:first-few} shows that the diameter of $\cK_n$ corresponds to the optimal path between one of the ``first few'' slow vertices. The last three sections use these ingredients to show both distributional convergence as well as the convergence of the moments of $\diam_w(\cK_n) -3\log{n}$ to the limiting random object thus completing the proof of the main result.

\subsection{Explicit distributions for distances between sets of vertices}
\label{sec:explicit}
In this section, we explain the proof by Janson of \eqref{eqn:sec-or-cij}.
We also extend that analysis to the smallest-weight path between disjoint sets of vertices. We remind the reader that the standing assumption henceforth is that each edge has exponential mean $n$ distribution. We start with the following lemma: 

\begin{lem}[Distances between sets of vertices]
\label{lem-dist-sets-vert}
Consider two disjoint non-empty sets $A,B\subseteq [n]$. Then,
	\eqn{
	d_w(A,B)\stackrel{d}{=} \sum_{k=|A|}^{N +|A| - 1} \frac{E_k}{k(n-k)},
	}
where 
\begin{itemize}
\item[(i)]$(E_k)_{k\geq 1}$ are i.i.d.\ mean $n$ exponential random variables;

\item[(ii)] $N$ is independent of the sequence $(E_k)_{k\geq 1}$ with the 
same distribution as the number of draws required to select the first black ball 
in an urn containing $|B|$ black balls and $n-|A|-|B|$ white balls, where one is drawing balls without replacement from the urn.
\end{itemize}
% $N=|A|+\min_{i=1}^{|B|} N_i$ where $(N_1,\ldots,N_{|B|})$ is a uniform choice of $|B|$ distinct elements from $[n-|A|]:= \set{1, 2, \ldots, n-|A|}$,  independent of the sequence $(E_k)_{k\geq 1}$.
\end{lem}

\proof We start exploring the neighborhood of the set $A$ in a similar way as in \cite{janson-cm-gp}. Recall that each edge has an exponential mean $n$ edge length. 
After having found the $\ell$th minimal edge and with $k=(|A|+\ell)$, 
there are $k(n-k)$ edges incident to the found vertices.
The minimal edge weight thus has an exponential distribution with mean  $n/k(n-k)$.
This process is stopped at the first time when we find a vertex in $B$. Since every new vertex added to the cluster of reached vertices is chosen uniformly amongst the set of present unreached vertices, the distribution of the number of steps required to reach a vertex in $B$ has the distribution 
$N$ asserted in the lemma, independently of the inter-arrival times of new vertices found. Thus the time it takes to find the first element in $B$ is 
	\eqn{
	\sum_{\ell=0}^{N-1} \frac{E_k}{(\ell+|A|)(n-\ell-|A|)}.
	}
Defining $k=\ell+|A|$ proves the claim.
\qed
\ \\ 

Now we specialize to a particular case of the above lemma.  Fix a vertex, \ch{say vertex $v=1$,} and another set $B\subseteq [n]\setminus \set{1}$. For much of the sequel we will be concerned with the optimal path between such a vertex and a set of size $|B| = \Theta(\sqrt{n})$. This is an appropriate time to think about two different but equivalent ways to find such an optimal path:
\ \\ 

\noindent{\bf Process 1: }
The first \ch{way to find the optimal path} is the exploration process described in the previous lemma  where we start at vertex $v=1$ and keep adding the closest vertex to the cluster until we hit a vertex in $B$. Write $M_B$ for the number of vertices other than $B$ \ch{that are found in this exploration}.  The previous lemma implies that
	\begin{equation}
	\label{eqn:exp-distrn-rel}
	\left(d_w(\set{1}, B), M_B\right) \stackrel{d}{=} \left(\sum_{k=1}^{N_B} \frac{E_k}{k(n-k)}, N_B\right),
	\end{equation} 
where $N_B$ is independent of the sequence $\ch{(E_k)_{k\geq 1}}$ and has the same distribution as the number of balls required to get the first black ball when drawing balls without replacement from an urn containing $|B|$ black balls and $n-1-|B|$ white balls.\ \\ 

\noindent{\bf Process 2:} The second \ch{way to find the optimal path} is the following. We think of water starting at source vertex $v=1$ at time $t=0$ percolating through the network at rate one using the edge lengths. Write $\SWG_t^{\sss(1)}$ (an acronym for the Smallest-Weight Graph) for the set of vertices reached by time $t$ starting from vertex $1$. More precisely,
	\begin{equation}
	\label{eqn:swg-def}
	\SWG_t^{\sss(1)}:= \set{u\in [n]\colon d_w(1,u)\leq t}.
	\end{equation} 
By \ch{convention}, vertex $v=1$ is in $\SWG_t^{\sss(1)}$ for all $t\geq 0$. Now note that the size process $(|\SWG_t^{\sss(1)}|)_{t\geq 0}$ is a pure-birth Markov process (with respect to the filtration $(\cF_t)_{t\geq 0} = (\sigma(\SWG_t))_{t\geq 0}$) with rate of birth given by $n/k(n-k)$ when the size $|\SWG_t^{\sss(1)}| =k$. Each new vertex added to this cluster is chosen uniformly amongst all available unreached vertices at that time, i.e. the vertices $[n]\setminus \SWG_t^{\sss(1)}$. Finally, the distance $d_w(\set{1}, B)$ can be recovered as 
	\begin{equation}
	\label{eqn:dwb-swg}
	d_w(\set{1}, B):= \inf\set{t\geq 0\colon \SWG_t^{\sss(1)}\cap B \neq \varnothing}.
	\end{equation}

In this section, we use Process 1 to prove the following initial result. We use Process 2 in Section \ref{sec:jt-con} below.
 
\begin{lem}[Distances between vertex and set of size $b\sqrt{n}$]
\label{lem-joint}
Let $B\subseteq [n]$ with $|B|=b\sqrt{n}$. Then as $n\to\infty$,
	\eqn{
	\left(d_w(\{1\},B)-\tfrac{1}{2}\log{n}, M_B/\sqrt{n}\right) \weakc \left(\Lambda+\log{(\hat{E}/b)},\hat{E}/b\right),
	}
where $\hat{E}$ is exponential with parameter 1, $\Lambda$ is Gumbel and $\hat{E}$ and $\Lambda$ are independent.
\end{lem}

\proof The above is equivalent to showing 
	\[(
	d_w(\set{1},B) - \log{M_B}, M_B/\sqrt{n} ) \weakc (\Lambda, \hat{E}/b),
	\]
with $\Lambda, \hat{E}$ independent standard Gumbel and $\exp(1)$ respectively. Fix constants $0< \alpha < \beta $ and $y\in \bR$. Define the event 
	\[
	A_n(y,\alpha, \beta):= \set{d_w(\set{1},B) - \log{M_B}\leq y}
	\cap\set{ \alpha \leq M_B/\sqrt{n} \leq \beta}.
	\]
Let $(E_k^\prime)_{k\geq 1}$ be independent sequence of {\bf mean one} exponential random variables.  Equation \eqref{eqn:exp-distrn-rel} implies 
	\begin{equation}
	\label{eqn:mult-p}
	\pr(A_n(y,\alpha, \beta)) 
	= \sum_{j=\alpha \sqrt{n}}^{\beta \sqrt{n}} 
	\pr\left(\sum_{k=1}^{j} \frac{n E_k^\prime}{k(n-k)} - \log{j}\leq y\right) \pr(N_B = j).
	\end{equation} 
Noting that  $\sum_{k=1}^j 1/j \approx \log{j} +\gamma$ as $j\to\infty$,  where $\gamma$ is Euler's constant, gives 
	\begin{equation}
	\label{eqn:sum-error-gumbel}
	\sum_{k=1}^{j} \frac{n E_k^\prime}{k(n-k)} - \log{j} \approx \sum_{k=1}^{j}\frac{E_k^\prime-1}{k} +\gamma +R_n,
	\end{equation} 
where the error term $R_n$ is independent of $j$ and is bounded by 
	\begin{equation}
	\label{eqn:rn-error-bd}
	|R_n| \leq \sum_{k=1}^{\beta \sqrt{n}} \frac{E_k^\prime}{n-k}\Probc 0,
	\end{equation} 
as $n\to\infty$. 
Thus, uniformly for $j\in [\alpha \sqrt{n}, \beta \sqrt{n}]$ 
	\[
	\pr\left(\sum_{k=1}^{j} \frac{n E_k^\prime}{k(n-k)} - \log{j}\leq y\right) 
	\to \pr\left(\sum_{k=1}^{\infty}\frac{E_k^\prime-1}{k} + \gamma \leq y \right).
	\]
It is easy to check (see e.g.\ \cite[Section 3]{janson-cm-gp}) that 
	\begin{equation}
	\label{eqn:gumbel-inf-exp}
	\sum_{k=1}^{\infty}\frac{E_k^\prime-1}{k} + \gamma \stackrel{d}{=} \Lambda.
	\end{equation}
By \eqref{eqn:mult-p} to complete the proof, it is enough to show that 
	\[
	\pr(\alpha \leq N_B/\sqrt{n}\leq \beta) \to \pr(\alpha \leq  \hat{E}/b\leq \beta).
	\]
This follows easily since for any $x > 0$
	\[
	\pr(N_B > x\sqrt{n}) = \prod_{k=1}^{x\sqrt{n}} \left(1-\frac{b\sqrt{n}}{n-1-k}\right) \sim \e^{- bx}, 
	\]
as $n\to\infty$. 
\qed

\subsection{{Poisson limit for the number of vertices with large minimal edge weights.}}
\label{sec:pp-limit}
The aim of this section is to understand the distribution of edges emanating from the slow vertices, namely the set of vertices for which the closest vertex is at distance $\approx \log{n}$.  For vertex $i\in [n]$, let $X_{\sss(i)}=\min_{j\in [n]} E_{ij}$ denote the minimal edge weight emanating from a given vertex $i\in [n]$.  Fix $\alpha \in \bR$ and let 
$N_n(\alpha)=\#\{i\in [n]\colon X_{\sss(i)}\geq \log{n}-\alpha\}$ denote the number of vertices 
with minimal outgoing edge weight at least $\log{n}-\alpha$. We prove the following Poisson limit for
$N_n(\alpha)$:

\begin{prop}[Number of vertices with large minimal edge weight]
\label{prop-large-min-edge-weight}
As $n\to\infty$,
	\eqn{
	N_n(\alpha)\weakc N(\alpha),
	}
where $N(\alpha)$ is a Poisson random variable with mean $\e^{\alpha}$. More precisely, 
	\begin{equation}
	\label{eqn:dtv}
	d_{\sss \mathrm{TV}}(N_n(\alpha), N(\alpha)) \leq \frac{2(1+\eps_n) \e^{2\alpha}\log{n}}{n},
	\end{equation}
where $d_{\sss \mathrm{TV}}$ denotes the total variation distance and $\eps_n = \exp\big(\frac{\log{n}-\alpha}{n}\big) - 1$.
\end{prop}

\proof We use the Stein-Chen method for Poisson approximation. Write 
	\[
	N_n(\alpha) = \sum_{i\in[n]} Z_i, \qquad  Z_i = \ind\set{X_{\sss(i)}\geq \log{n}-\alpha}. 
	\] 
For fixed $i\in [n]$, note that $X_{\sss(i)}$ has an exponential distribution with mean $n/(n-1)$. Writing $p_n = \pr(Z_i=1)$ so that $\lambda:=\E(N_n(\alpha)) = np_n$, it is easy to check that 
	\begin{equation}
	\label{eqn:lambda-comp}
	\E(N_n(\alpha)) = (1+\eps_n) \e^{\alpha}. 
	\end{equation}
Thus, $\lambda \to \e^{\alpha}$ as $n\to\infty$. For each fixed $i\in [n]$, suppose we can couple $N_n(\alpha)$ with a random variable $W^\prime_i$ such that the marginal distribution of $W^\prime_i$ is 
	\begin{equation}
	\label{eqn:wprime-dist}
	W^\prime_i+1 \stackrel{d}{=} N_n(\alpha)\bigl|_{\set{Z_i =1}},
	\end{equation}
i.e.,  $W^\prime_i+1$ has the same distribution as $ N_n(\alpha)$ conditionally on $\set{Z_i =1}$.
Then Stein-Chen theory \cite{MR1163825} implies that in total variation distance 
	\begin{equation}
	\label{eqn:dtv-bound}
	d_{\sss \mathrm{TV}}(\cL(N_n(\alpha)),\mathrm{Poi}(\lambda)) 
	\leq (1\wedge\lambda^{-1}) \sum_{i\in[n]} \E(Z_i) \E(|N_n(\alpha)-W^\prime_i|)
	\end{equation} 
Let us describe $W_1^\prime$, the same construction switching indices works for any $i$. Let $\sS_n:=\set{\cK_n, (E_e)_{e\in \cE_n}}$ be the original edge lengths and let $N_n(\alpha)$ be defined as above for the random metric space $\sS_n$.   Let us construct the edge lengths of $\cK_n$ conditional on the event $\set{Z_1=1}$ so that $X_{\sss(1)}-\log{n}\geq -\alpha$. We shall write $\sS_n^\prime:=\set{\cK_n, (E_e^\prime)_{e\in \cE_n}}$ for $\sS_n$ conditioned on this event.  Note that this event only affects edges incident to vertex $1$ and further, by the lack of memory property of the exponential distribution, every such edge incident to vertex $1$ has distribution $\log{n}-\alpha +E$ where $E$ is an exponential mean $n$ random variable, independently across edges. Thus, we can construct the edge lengths on $\sS_n^\prime$ using the edge lengths $E_e$ in $\sS_n$ by the following description: 
	\begin{enumeratea}
	\item For each edge $e=\set{1,i}$ incident to vertex $i$, set $E_e^\prime = \log{n}-\alpha +E_e $.
	\item For any edge not incident to vertex $1$, set $E_e^\prime = E_e$. 
	\end{enumeratea} 
Define $X_{\sss(i)}^\prime$ analogously to $X_{\sss(i)}$ as the minimal edge length incident to vertex $i$ but in $\sS_n^\prime$. Finally, define 
	\[
	Z_i^\prime:= \ind\set{X_{\sss(i)}^\prime> \log{n} - \alpha}, \qquad 
	W_1^\prime = \sum_{v\neq 1} \ind\set{X_{\sss(v)}^\prime\geq \log{n}-\alpha}. 
	\] 
Then $W^\prime_1$ by construction has the required distribution in \eqref{eqn:wprime-dist}. Note that
	\[
	|N_n(\alpha) - W_1^\prime|\leq  
	\ind\set{X_{\sss(1)}> \log{n} - \alpha} + \sum_{i\neq 1}|Z_i-Z_i^\prime|.
	\]
Taking expectations, by symmetry,
	\begin{equation}
	\label{eqn:sum-bd}
	\E(|N_n(\alpha)-W^\prime_1|)\leq p_n + (n-1) \E|Z_{2} - Z_{2}^{\prime}|.
	\end{equation}
Now 
	\[
	\E|Z_{2} - Z_{2}^{\prime}| = \pr(Z_{2}=1,~ Z_{2}^{\prime}=0)+ \pr(Z_{2} =0,~ Z_{2}^{\prime}=1).
	\] 
Since the edge lengths in $\sS_n^\prime$ are at least as large as the edge lengths in $\sS_n$, we have $\set{Z_{2}=1,~ Z_{2}^{\prime}=0} = \varnothing $. For the second term 
	\[
	\set{Z_{2} =0,~ Z_{2}^{\prime}=1}\equiv \set{E_{2,1}< \log{n}-\alpha, 
	\min_{j\neq 1,2} E_{2,j} \geq \log{n}-\alpha}.
	\] 
Since $E_{i,j}$ are exponential mean $n$, we immediately get
	\[
	\pr(Z_{2} =0,~ Z_{2}^{\prime}=1) \leq \frac{\e^{\alpha}\log{n}}{n^2}.
	\]
Using this in \eqref{eqn:sum-bd}, the total variation bound \eqref{eqn:dtv-bound} completes the proof. \qed
	% \eqn{
	% N_n(\alpha)=\sum_{i\in [n]} \indic{X_{\sss(i)}\geq \log{n}-\alpha}.
	% }
	
%\todo{Update from here!}

\subsection{{Joint convergence of distances between multiple vertices}}
\label{sec:jt-con}
The aim of this section is to understand the re-centered asymptotic joint distribution of the \ch{minimal weight} between multiple vertices. To prove this, it turns our that Process 2 using the \ch{smallest-weight graph} $\SWG_t^{\sss(v)}$ from vertices $v\in [n]$ is more useful than Process 1. Versions of Proposition \ref{prop-many-vertices-weights} below has appeared before in \cite{aldous-exch,shankar-er,david-shankar}. We give a new proof, \ch{both for completeness as well as since we} need a variant of this argument in the sequel. 

Fix $m\geq 2$. Let $(\Lambda_\alpha)_{\alpha\in [m]}$ and $(\Lambda_{\alpha \beta})_{\alpha,\beta\in [m],s<t}$ 
be independent standard Gumbel random variables. In the following proposition, we identify the limiting distribution of $(d_w(\alpha,\beta)-\log{n})_{\alpha,\beta\in [m],\alpha<\beta}$, an extension of the result given in \eqref{eqn:sec-or-cij} proved by Janson \cite{janson-cm-gp} for $m=2$:

\begin{prop}[Joint distances between many vertices]
\label{prop-many-vertices-weights}
As $n\to\infty$,
	\eqn{
	(d_w(\alpha,\beta)-\log{n})_{\alpha,\beta\in [m],\alpha<\beta}
	\weakc (\Lambda_\alpha+\Lambda_\beta-\Lambda_{\alpha \beta})_{\alpha,\beta\in [m],\alpha<\beta}.
	}
\end{prop}

\proof Fix $m\geq 2$.  Write 
	\begin{equation}
	\label{eqn:mvd-def}
	\mvD(m):= (\Lambda_\alpha+\Lambda_\beta-\Lambda_{\alpha \beta})_{\alpha,\beta\in [m],\alpha<\beta},
	\end{equation} 
for the limiting array. The idea of the proof is as follows. We start by sequentially growing the smallest-weight graphs $\SWG$'s from the $m$ vertices until they meet. This gives us a sequence of collision times $(T_{\alpha \beta})_{\alpha < \beta \in [m]}$.  An appropriately chosen linear transformation of these collision times stochastically dominates the array of the lengths of shortest paths. We show that this linear transformation of the collision times converges to the array $\mvD$.   A simple limiting argument using the convergence of the marginal distribution of two point distances implies that the joint distribution of the distances themselves converge to $\mvD$ and this completes the proof.

Let us now start with the proof. Throughout we write $\sS_n$ for the random metric space $(\cK_n, \set{E_e}_{e\in \cE_n})$, where once again we remind the reader that $E_e$ are i.i.d.\ exponential random variables with mean $n$.  Now start the smallest weight cluster $\SWG^{\sss(1)}_t$ from vertex $\alpha = 1$.  Write 
	\eqn{
	T_1=\inf\{t\colon |\SWG^{\sss(1)}_t|=\sqrt{n}\}
	}
\ch{for the time for $\SWG^{\sss(1)}_t$ to grow to size $\sqrt{n}$.}
Then, since $T_1\stackrel{d}{=} \sum_{k=1}^{\sqrt{n}} nE_k/[n(n-k)]$, this implies (see \eqref{eqn:sum-error-gumbel} and \eqref{eqn:gumbel-inf-exp}) that 
	\eqn{
	T_1-\tfrac{1}{2}\log{n}\weakc \log(1/\hat{E}_1),
	\label{eqn:e1-def}
	}
where $\hat{E}_1$ is exponential with mean 1. For every vertex $v\in \SWG_t^{\sss(1)}$, write $B^{\sss(1)}(v) := d_w(1,v)$ for the time when the flow from vertex $1$ reaches $v$. \ch{We now work conditionally} on the flow cluster $\SWG^{\sss(1)}_{T_1}$. By construction, as $n\to\infty$,
 	\begin{equation}
 	\label{eqn:2-not-in-1}
 	\prob(2\notin \SWG_{T_1}^{\sss(1)})= 1-\frac{\sqrt{n}}{n} \to 1.
 	\end{equation}
Further, by the \ch{memoryless} property of the exponential distribution, \ch{conditionally} on $\SWG_{T_1}^{\sss(1)}$, for every boundary edge $e= \set{u,v}$ with $u\in \SWG_{T_1}^{\sss(1)}$ and $v\notin \SWG_{T_1}^{\sss(1)}$, the remaining edge length $E_e - (T_1 - B^{\sss(1)}(u))$ has an exponential distribution with mean $n$, \ch{and all these remaining edge lengths are independent.}

Freeze the cluster $\SWG_{T_1}^{\sss(1)}$. Start a flow from vertex $2$ as the source and write $\SWG_t^{\sss(2)}$ for the \ch{smallest-weight} graph. Write 
	\begin{equation}
	\label{eqn:t12-def}
	T_{12}:=\inf\set{t: \SWG_t^{\sss(2)} \cap \SWG_{T_1}^{\sss(1)} \neq \varnothing},
	\end{equation}
\ch{so that $T_{12}$ is} the first time that a vertex in the flow cluster from vertex $\alpha=1$ at time $T_1$ is hit by the flow cluster from $2$. \ch{Conditionally} on $\SWG_{T_1}^{\sss(1)}$, on the event $\set{2\notin \SWG_{T_1}^{\sss(1)} }$ we have that
\begin{enumeratea}
	\item the \ch{smallest-weight} path between $1$ and $2$ is given by $d_w(1,2) = T_1 + T_{12}$. 
	\item the random variable $T_{12}$ has the same distribution as $d_w(\set{1}, B)$ in the random (unconditional)  metric space $\sS_n$ where $B$ is a fixed set of size $\sqrt{n}$.  
\end{enumeratea}

By Lemma \ref{lem-joint} with $b=1$ we immediately get 
	\eqn{
	(T_{12}-\tfrac{1}{2}\log{n}, |\SWG_{T_{12}}^{\sss(2)}|/\sqrt{n})\weakc (\log{(1/\hat{E}_2)}+\log{(\hat{E}_{12})},\hat{E}_{12}),
	\label{eqn:e2-e12-def}
	}
where $\hat{E}_2$ and $\hat{E}_{12}$ are independent of $\hat{E}_1$ in \eqref{eqn:e1-def}.  Combining \eqref{eqn:e1-def} and \eqref{eqn:e2-e12-def} we get 
	\eqan{
	(d_w(1,2)-\log{n}, |\SWG_{T_{12}}^{\sss(2)}|/\sqrt{n})&
	=(T_{12}-\tfrac{1}{2}\log{n}+T_1-\tfrac{1}{2}\log{n}, N/\sqrt{n}) \notag \\
	&\weakc (\log(1/\hat{E}_1)+\log(1/\hat{E}_2)+\log(\hat{E}_{12}), \hat{E}_{12}).\label{eqn:convg-m-2}
	}
This proves the claim for $m=2$. We next extend the computation to $m=3$. 

For ease of notation, write $\cB = \sqrt{n} = |\SWG_{T_1}^{\sss(1)}|$ and 
$\cR = |\SWG_{T_{12}}^{\sss(2)}|$, here $\cB$ and $\cR$ will be mnemonics for ``black'' and ``red'' respectively. We now work \ch{conditionally} on $\cA := \SWG_{T_1}^{\sss(1)} \cup \SWG_{T_{12}}^{\sss(2)}$. Since $|\cA| = \Theta_P(\sqrt{n})$, 
	\begin{equation}
	\label{eqn:3-not-in-12}
	\prob(3\notin \SWG_{T_1}^{\sss(1)} \cup \SWG_{T_{12}}^{\sss(2)} )\to 1 \qquad \mbox{as } n\to\infty.
	\end{equation}
Freeze the above two flow clusters. Start a flow from vertex $\beta =3$ and consider the \ch{smallest-weight graph} $\SWG_{t}^{\sss(3)}$ emanating from vertex $3$. We need to modify this process after the first time it finds a vertex in $\cA= \SWG_{T_1}^{\sss(1)} \cup \SWG_{T_{12}}^{\sss(2)}$, namely after time 
	\[
	T_3^* = \inf\set{t: \SWG_t^{\sss(3)} \cap \cA \neq \varnothing}.
	\]
Suppose this happens \ch{due to} $\SWG_{T_3}^{\sss(3)}$ finding a vertex in $\SWG_{T_1}^{\sss(1)}$. Remove all vertices in $\SWG_{T_1}^{\sss(1)}$ and all adjacent edges from $\cK_n$ and then continue until the process finds a vertex in $\SWG_{T_{12}}^{\sss(2)}$. Similarly if this happens \ch{due} to a vertex in $\SWG_{T_{12}}^{\sss(2)}$ being found, then remove all vertices in  $\SWG_{T_{12}}^{\sss(2)}$ and continue. Although this is not quite the \ch{smallest-weight graph} emanating from vertex $3$, to minimize notational overhead, we shall continue to denote this modified process by the same $\set{\SWG_t^{\sss(3)}}_{t\geq 0}$. Define the stopping times
	\[
	T_{13} = \inf\set{t\geq 0: \SWG_t^{\sss(3)}\cap \SWG_{T_1}{\sss(1)} \neq \ch{\varnothing}},
	\]
and 
	\[
	T_{23} = \inf\set{t\geq 0: \SWG_t^{\sss(3)}\cap \SWG_{T_{12}}{\sss(2)} \neq \ch{\varnothing}}.
	\]
Similarly, define the sizes of the cluster $\SWG_t^{\sss(3)}$ at these stopping times as
	\eqn{
	C_n^{\sss(13)}=|\SWG_{T_{13}}^{\sss(3)}|, \qquad 	C_n^{\sss(23)}=|\SWG_{T_{13}}^{\sss(3)}|.
	}
Similar to the urn description in \eqref{eqn:exp-distrn-rel}, it is easy to check that conditionally on $\cA$ \ch{and} on the event $\set{3\notin \cA}$, the distribution of the random variables $(T_{13}, T_{23}, C_n^{\sss(13)}, C_n^{\sss(23)})$ can be constructed as follows:\\
Consider an urn with $n$ balls out of which $\cB = |\SWG_{T_1}^{\sss(1)}|$ black balls, $\cR = |\SWG_{T_{12}}^{\sss(2)}|$ red balls and the remaining $n-\cB-\cR$ white balls. Also let $(E_k)_{k\geq 1}$ be an independent sequence of mean $n$ exponential random variables.  Start drawing balls at random without replacement till the first time $\cN_1$ that we get either a black or a red ball. 
	\begin{enumeratea}
	\item Suppose the first ball amongst the black or red balls is a black ball. Remove all black balls so that there are now
	 $(n-\cN_1 - \cB)$ balls in the urn. Continue drawing balls without replacement till we get a red ball. 
	Let $\cN_2> \cN_1$ be the time for the first pick of a red ball. Let $C_n^{\sss(13)} = \cN_1$, $C_n^{\sss(23)} = \cN_2$.
 	Finally, let 
		\begin{equation}
		T_{12}:= \sum_{k=1}^{\cN_1} \frac{ E_k}{k(n-k)}, \qquad
		T_{23}:= T_{12} + \sum_{k=\cN_1+1}^{\cN_2} \frac{ E_k}{k(n- k -\cB)}
		\end{equation} 
		where as before, $(E_k)_{k\geq 1}$ is an independent sequence of exponential random variables with mean $n$. 
	\item Suppose the first ball amongst black and red balls to be picked is a red ball. 
	Then, in the above formulae, \ch{simply} interchange the roles of $1$ and $2$ and $\cB$ and $\cR$.  
	\end{enumeratea} 
Using \eqref{eqn:e2-e12-def} and arguing exactly as in the proof of Lemma \ref{lem-joint}, \ch{we see} that 	
	
	\begin{align}
		\label{eqn:e123-t-convg}
		\biggl(\frac{C_n^{\sss(13)}}{\sqrt{n}}, \frac{C_n^{\sss(23)}}{\sqrt{n}},
		&T_{13}-\tfrac{1}{2}\log{n}, T_{23}-\tfrac{1}{2}\log{n}\biggr)\weakc\\ 
		& (\hat{E}_{13}, \hat{E}_{23}/\hat{E}_{12}, \log(1/\hat{E}_3)+\log(\hat{E}_{13}), \log(1/\hat{E}_3)+\log(\hat{E}_{23}/\hat{E}_{12} ) \notag
	\end{align}
Here $\hat{E}_3, \hat{E}_{13}, \hat{E}_{23}$ are independent of $\hat{E}_1, \hat{E}_2, \hat{E}_{12}$ and i.i.d.\ exponential mean-one random variables. Now note that by construction, there is a path of length $D_n(1,3):= T_1+ T_{13}$ between vertices $1$ and $3$ and similarly of length $D_n(2,3): = T_{12}+ T_{23}$ between vertices $2$ and $3$.   
% Define
% 	\eqn{
% 	T_{13}=T^{\sss(3)}_{C_n^{\sss(13)}},
% 	\qquad
% 	T_{23}=T^{\sss(3)}_{C_n^{\sss(23)}}.
% 	}
% Then, conditionally on $E_{12}$,
% 	\eqn{
% 	()
% 	\weakc (-\log(E_3)+\log(E_{13}), -\log(E_3)+\log(E_{23}/E_{12}).
% 	}
Thus, by \eqref{eqn:e2-e12-def} and \eqref{eqn:e123-t-convg}
	\eqn{
	d_w(1,3)-\log{n}\leq T_{13}-\tfrac{1}{2}\log{n}+T_1-\tfrac{1}{2}\log{n}
	\weakc \log(1/\hat{E}_1)+\log(1/\hat{E}_3)+\log(\hat{E}_{13}),
	}
and
	\eqan{
	d_w(2,3)-\log{n}& \leq T_{23}-\tfrac{1}{2}\log{n}+T_{12}-\tfrac{1}{2}\log{n}\\
	&\weakc \log(1/\hat{E}_3)+\log(\hat{E}_{23}/\hat{E}_{12}+\log(1/\hat{E}_2)+\log(\hat{E}_{12})\nn\\
	&=\log(1/\hat{E}_2)+\log(1/\hat{E}_3)+\log(\hat{E}_{23}),\nn
	}
Thus the limiting array $\mvD(3)$ in \eqref{eqn:mvd-def} is a limiting upper bound in the weak sense for the array $\mvd_n(3):= (d_w(\alpha,\beta) - \log{n}: 1\leq \alpha < \beta \leq 3)$. However, we have equality for $m=2$ by \eqref{eqn:convg-m-2}. Thus the marginals of $\mvd_n(3)$ converge to the marginals of $\mvD$ as $n\to\infty$. This implies $\mvd_n(3) \weakc \mvD(3)$ as $n\to\infty$.  

This entire construction extends inductively for higher values of $m$ and thus completes the proof. 
\qed

\noindent {\bf Remark.} We learned about this reduction from the sums of collision times to lengths of optimal paths via stochastic domination from \cite{salezjoint}. 
\ \\ 

The following is an easy corollary of the proof of the above result. Recall that for any 2 vertices $\alpha, \beta \in [n]$, $\mvpi(\alpha, \beta)$ denotes the unique shortest path (geodesic) between them. 
\begin{cor}
\label{cor:mult-path-prop}
Consider the random metric space $\sS_n = (\cK_n, \set{E_e}_{e\in \cE_n})$. Fix $m\geq 2$. Then,
\begin{enumeratea}
	\item Let $D_n$ be the event that $\exists \alpha\neq \beta\neq \gamma \in [m]$ such that $\gamma \in \mvpi(\alpha, \beta)$. Then $\pr(D_n)\to 0$ as $n\to\infty$. 
	\item Fix $1/2< \vartheta <1$. Consider the \ch{smallest-weight graph}s $\set{\SWG_{\vartheta\log{n}}^{\sss(i)}}_{i\in[m]}$ from these $m$ vertices \ch{at time} $\vartheta \log{n}$. % Consider the subset $\cA([m],\alpha\log{n} ) = \cup_{i=1}^m \SWG_{\alpha \log{n}}^{\sss(i)} $.
	 Then whp, the shortest paths $\mvpi(\alpha,\beta)$ are contained in the union of these balls, \ch{i.e.,} as $n\to\infty$,
	\[
	\pr(\mvpi(\alpha, \beta) \subseteq \cup_{i=1}^m \SWG_{\vartheta \log{n}}^{\sss(i)} ~ \forall \alpha, \beta \in [m]) \to 1.
	\] 
\end{enumeratea}
\end{cor}
\proof Part(a) follows from extending \eqref{eqn:2-not-in-1} and \eqref{eqn:3-not-in-12} to general $m$. Part (b) follows from the above proof which proves that for any pair of vertices $\alpha, \beta$, $\pi(\alpha, \beta)$ can be found in $\SWG_{r_n}^{\sss(\alpha)}\cup \SWG_{r_n}^{\sss(\beta)} $ where $r_n  = \frac{1}{2}\log{n} + O_P(1)$. 
\qed

\subsection{{Distances between vertices with large minimal edge weight}}
\label{sec:dis-min-mult}
Fix $\alpha \in \bR$. Recall that $N_n(\alpha) = \sum_{i=1}^n \ind\set{X_{\sss(i)}\geq \log{n} - \alpha}$ denotes the number of vertices with minimum outgoing edge length at least $\log{n} - \alpha$. Fix $m\geq 2$ and condition on the event $N_n(\alpha)=m$. Let $V_1, \ldots, V_m$ denote the $m$ vertices for which 
$X_{\sss(V_i)}\geq \log{n}-\alpha$.

\ch{Our aim in this section} is to understand, \ch{conditionally} on the event $\set{N_n(\alpha) =m}$, the asymptotic joint distribution of $(d_w(V_i, V_j): i< j \in [m])$. Recall the array $\mvD(m)$ from \eqref{eqn:mvd-def} giving the asymptotic joint distribution of the re-centered (by $\log{n}$) length of smallest paths between $m$ typical vertices in $\sS_n$. 
The main aim of this section is to prove the following result:  
 % We next extend Proposition \ref{prop-many-vertices-weights}
 % to this conditional setting:

\begin{prop}[Distances between vertices with large minimal edge weight]
\label{prop-dist-large-min-edge-weight}
Fix $\alpha\in \bR$ and $m\geq 2$.  Conditionally on $N_n(\alpha)=m$, as $n\to\infty$,
	\eqn{
	(d_w(V_i,V_j)-3\log{n}+2\alpha)_{i,j\in [m],i<j}
	\weakc (\Lambda_i+\Lambda_j-\Lambda_{ij})_{i,j\in [m],i<j}:= \mvD(m).
	}
\end{prop}

%\todo{Update from here!}
\proof Let us \ch{start by disentangling} exactly what the conditioning event $\set{N_n(\alpha) = m}$ implies about the edge length distribution. We write $\sS_n^\prime(\tra, \tca)$ for the conditioned metric space. Here ``$\tra, \tca$'' are short for ``translation'' and ``conditioning'' respectively. This will become clear below.  The basic idea is to use our original (unconditioned) random metric space $\sS_n$ to generate the metric space $\sS_n^{\prime}(\tra, \tca)$.  To ease notation, we assume w.l.o.g.\ that $V_i = i$. Then this conditioning implies that the edge lengths of $\sS_n^\prime(\tra, \tca)$ can be constructed \ch{by} the following two rules:  
	\begin{enumeratea}
	\item {\bf Translation:} Every edge $E_e^\prime$ incident to one of the vertices in $[m]$ is conditioned to be at least $\log{n} - \alpha$. By the \ch{memoryless property of the exponential distribution,} we can write $E_e^\prime = \log{n} - \alpha + E_e$ where $\ch{(E_e)}$ are \ch{an independent} family of mean $n$ independent exponential random variables. 
	\item {\bf Conditioning:} For every vertex $i\notin [n]\setminus [m]$, the edges $\ch{(E_{i,j}^\prime)_{ j\notin [m]}}$ are independent exponential mean $n$ random variables conditioned on 
	\begin{equation}
	\label{eqn:cond}
		X_{\sss(i), [m+1:n]}: = \ch{\min_{m+1\leq  j\leq n} E_{i,j}^\prime} <  \log{n}- \alpha.
	\end{equation} 
\end{enumeratea}
Let us use our original metric space $\sS_n$ to sequentially overlay the effect of the above 2 events. More precisely, we will use our original metric space $\sS_n$ to construct $\sS_n^\prime(\tra,\tca)$ in two steps.  Recall that we \ch{have} used $\mvpi(i, j)$ for the \ch{smallest-weight} path between $i, j$ in $\sS_n$. \csb{The following lemma deals with the effect of the simpler translation event (without dealing with the conditioning),} \crh{and will be the starting point of our analysis:}

\begin{lem}
\label{lem:translation}
Fix $m\geq 1$ and consider the metric space $\sS_n$. For every edge $e$ incident to one of the vertices \ch{in} $[m]$, replace the edge $E_e$ by $E_e + \log{n} - \alpha$. Leave all other edges unchanged. Call this new metric space $\sS_n^\prime (\tra)$. Write $\mvpi^\prime(i, j)$ for the \ch{smallest-weight} path between $i, j $ and write $d_w^\prime$ for the corresponding metric. Then, for all $i, j \in [m]$, 
	\begin{equation}
	\label{eqn:path-relt}
	\mvpi^\prime(i, j) = \mvpi(i, j), \qquad  d_w^\prime(i, j) = d_w(i, j) + 2\log{n}-2\alpha. 
	\end{equation}
 In particular,
	\[\left(d_w^\prime(i,j) -3\log{n}+2\alpha\right)_{i,j\in [m],i<j}
	\weakc \left(\Lambda_i+\Lambda_j-\Lambda_{ij}\right)_{i,j\in [m],i<j}. 
	\]
\end{lem}
\proof The distributional convergence follows from \eqref{eqn:path-relt} and Proposition \ref{prop-many-vertices-weights}. Equation \eqref{eqn:path-relt} follows since we can construct the \ch{smallest-weight} path problem for $\sS_n^\prime(\tra)$ as follows. To $\sS_n$ adjoin $m$ new vertices $\set{i^\prime: i'\in [m]}$. Each new vertex $i^\prime$ has only one edge, \ch{namely,} to vertex $i$ of length $\log{n} - \alpha$. Call this new metric space $\sS_n^*$ and the corresponding metric $d_w^*$ and \ch{smallest-weight} path $\mvpi^*(\cdot, \cdot)$. \csb{Then the metric space $\sS_n^\prime(\tra)$ can be constructed as follows:} For $i,j\in [m]$ let $d_w^\prime(i,j) = d_w^*(i^\prime, j^\prime)$ and $\mvpi^*(i^\prime, j^\prime) = \set{i^\prime \leadsto i} \cup \mvpi^\prime(i,j) \cup \set{j\leadsto j^\prime}$. 
\qed
\bigskip
% \todo{Remco: I do not get why Lemma \ref{lem:translation} is related to the original problem? We never refer to it either...}

Let us now construct the full metric space $\sS_n^\prime(\tra, \tca)$. We construct this from $\sS_n$ in 4 steps. Fix $1/2<\vartheta<1$. Write $\cB_n(\alpha) = \set{v\in [n]\ch{\setminus[m]}\colon X_{\sss(v)} \geq \log{n} -\alpha}$. 
This is the set of ``bad'' vertices whose edges we need to ``correct''. 

	\begin{enumeratea}
	\item First construct the \ch{smallest-weight graph}s $\ch{\big(\SWG_{\vartheta \log{n}}^{\sss(i)}\big)_{i\in[m]}}$. By Corollary \ref{cor:mult-path-prop}, with high probability $\mvpi(i,j) \subseteq \cup_{i=1}^m \SWG_{\vartheta \log{n}}^{\sss(i)} $ for all $i,j\in [m]$. 

	\item Now reveal all the other edges. 

	\item {\bf Translation:} To each edge incident to one of the vertices $i\in [m]$, add $\log{n}-\alpha$. \csb{This gives us the metric space $\sS_n'(\tra)$. The effect of this has been analyzed in Lemma \ref{lem:translation}.}

	\item {\bf Conditioning:}  Now consider the vertices \ch{in} $\cB_n(\alpha)$. Note that by \ch{Proposition \ref{prop-large-min-edge-weight} and} as $n\to\infty$, $|\cB_n(\alpha)|\weakc \mbox{Poi}(\e^{\alpha})$.  When  $\mvpi(i,j) \subseteq \cup_{i=1}^m \SWG_{\vartheta \log{n}}^{\sss(i)}$ for all $i,j\in [m]$ then 
		\[
		\cB_n(\alpha)\cap\cup_{i=1}^m \SWG_{\vartheta \log{n}}^{\sss(i)}= \varnothing,
		\]
	 since $\vartheta < 1$ and thus every vertex $v\in \cup_{i=1}^m \SWG_{\vartheta \log{n}}^{\sss(i)}$ 
	has at least one edge with length $\leq \vartheta \log{n}$. To complete the construction, \ch{we} 
	resample the edge lengths $\ch{(E_{v, i})_{v\in \cB_n(\alpha), m+1\leq i\leq n}}$ such that for every 
	vertex $v\in \cB_n(\alpha)$, we have $X_{\sss(v), [m+1: n]}< \log{n}-\alpha$. 
	\end{enumeratea} 
This completes the construction of $\sS_n^\prime(\tca, \tra)$. Now, after resampling, for $v\in \cB_n(\alpha)$ \ch{and} $i\geq m+1$, we write $E^\prime_{v, i}$ for the re-sampled edge lengths. For $v\in \cB_n\ch{(\alpha)}$, write 
	\[
	X^*_{\SWG}(v) = \min_{j\in  \cup_{i=1}^m \SWG_{\vartheta \log{n}}^{\sss(i)}} E^\prime_{v, j}
	\]
%\todo{Remco: Do you mean $E_{v, j}'$ here?}
for the \ch{smallest} edge weight from $v$ to $\cup_{i=1}^m\SWG_{\vartheta \log{n}}^{\sss(i)}$. We shall show that
	\begin{equation}
	\label{eqn:to-infty}
 	\frac{X^*_{\SWG}(v)}{\vartheta \log{n}}\to \infty \qquad \mbox{as } n\to\infty.
	\end{equation}
This implies that whp the resampling of the edge lengths of $v$ does not disturb $\cup_{i=1}^m\SWG_{\vartheta \log{n}}^{\sss(i)}$ and, in particular, the \ch{smallest-weight} path between $i,j$ in $\sS_n^\prime(\tra,\tca)$ for all $i,j\in [m]$ \csb{is the same as that in $\sS_n^\prime(\tra)$. Lemma \ref{lem:translation} then completes the proof.} 

We now show \eqref{eqn:to-infty}. Let us first estimate the size of $|\SWG_t^{\sss(i)}|$. Recall from Section \ref{sec:explicit} that, for any $t\geq 0$ and any $i\in [m]$,
	\[
	|\SWG_t^{\sss(i)}| \stackrel{d}{=} 1+ \max \set{l\geq 1: \sum_{k=1}^l \frac{n E_k^\prime}{k(n-k)} \leq t }.
	\]
Here $(E_k^\prime)_{k\geq 1}$ is an \ch{i.i.d.\ sequence of} exponential mean one \ch{random variables.} Obviously, this process is stochastically dominated by the process 
	\[
	Y(t):= 1+ \max\set{l\geq 1: \sum_{k=1}^l \frac{E_k^\prime}{k} \leq t }. 
	\]
The process $(Y(t))_{t\geq 1}$ is called the Yule process and is one 
of the standard examples of a pure birth process. In particular, (see \ch{e.g.} 
\cite{norris1998markov}), $(\e^{-t}Y(t))_{t\geq 0}$ is an $\bL^2$-bounded 
positive martingale. Therefore, for any $\vartheta^\prime > \vartheta$ as $n\to\infty$, 
	\[
	\frac{|\SWG^{\sss(i)}_{\vartheta \log{n}}|}{n^{\vartheta^\prime}}
	\leq \frac{Y(\vartheta \log{n})}{n^{\vartheta^\prime}} \Probc 0.
	\]
As a result, $|\cup_{i=1}^m\SWG_{\vartheta\log{n}}^{\sss(i)}| = o_P(n^{\vartheta^\prime})$. 
The following simple lemma which we give without proof, completes the proof of \eqref{eqn:to-infty} 
and thus the proof of  Proposition \ref{prop-dist-large-min-edge-weight}:

	\begin{lem}
	\label{lem:to-inf}
	Let $D_1, D_2, \ldots, D_n$ be \ch{i.i.d.\ exponential} mean $n$ random variables 
	conditioned on $X_{\sss(1)} = \min_{1\leq i\leq n} D_i < \log{n} - \alpha$. 
	Let $X^* = \min_{1\leq i\leq n^{\vartheta^\prime}} D_i$. Then, with $W\sim \exp(1)$,
		\[
		\frac{X^*}{n^{1-\vartheta^\prime}} \weakc W \qquad \mbox{ as } n\to\infty.
		\]
	\end{lem}\qed

\subsection{{Reduction to distances between vertices with large minimal edge weights}}
\label{sec:first-few}
The previous section analyzed distances between the vertices whose minimal outgoing edge is large (like $\log{n}+O_P(1)$). The distances between these vertices \ch{are then close to} $3\log{n}+O_P(1)$.  The aim of this section is to show that these are the only vertices that matter \ch{for the weight diameter. We achieve this} by considering distances between vertices whose minimal outgoing edge is ``small'' and \ch{showing} that the distance between such vertices are not large enough to create the diameter and thus can be ignored. 

\ch{We start with some notation.} Fix $\alpha >0$ and define
	\eqn{\label{eqn:rn-def}
	R_n(\alpha)
	=\#\{i,j\in [n]\colon X_{\sss(i)}\leq \log{n}-\alpha, X_{\sss(j)}\leq \log{n}+\alpha/2, d_w(i,j)\geq 3\log{n}-\alpha/8\}.
	}
\ch{The random variable $R_n(\alpha)$ counts} the number of ordered pairs of vertices $(i,j) \in [n]\times [n]$ \ch{that satisfy that} the minimal outgoing edge of vertex $i$ is less than $\log{n}-\alpha$, the minimal outgoing edge of $j$ is less than $\log{n}+\alpha/2$ and yet the distance between $i,j$ is greater than $\log{n} -\alpha/8$. 
\ch{The following lemma gives an upper bound on the expected value of $R_n(\alpha)$:}

\begin{prop}[Distances from vertices with small minimal weight]
\label{prop-small-min-edge-weight} \ \\ 
There exists a constant $C> 0$ such that for all $\alpha > 0$,
	\eqn{
\limsup_{n\to\infty}	\E[R_n(\alpha)]\leq C \e^{-\alpha/16}.
	}
\end{prop}

\proof We compute
	\eqn{
	\E[R_n(\alpha)]=n^2\pr(d_w(1,2)\geq 3\log{n}-\alpha/8, X_{\sss(1)}\leq \log{n}-\alpha, X_{\sss(2)}\leq \log{n}+\alpha/2).
	}
Note that 
	\[(X_{\sss(1)},X_{\sss (2)})\stackrel{d}{=}
	\left(\min \left[\frac{n}{n-2}E_1^*,nE_{12}^*\right], \min\left[\frac{n}{n-2}E_2^*, nE_{12}^*\right]\right),
	\] 
where $E_1^*, E_2^*, E_{12}^*$ are independent exponential random variables with mean 1. Here $nE_{12}^*$ represents the weight of the direct edge between vertices $1,2$, while \ch{for $i\in \{1,2\}$, $nE_i^*/(n-2)$ represents the minimal outgoing edges from vertex $i$ to the remaining vertices $[n]\setminus \set{1,2}$.} 
 
On the event $\set{d_w(1,2)\geq 3\log{n}-\alpha/8}$, we have that
$nE_{12}^*\geq d_w(1,2)\geq 3\log{n}-\alpha/8.$
As a result, when $d_w(1,2)\geq 3\log{n}-\alpha/8$, {\bf unless} 
	\begin{equation}
	\label{eqn:unless}
	\max( \frac{n}{n-2}E_1^*, \frac{n}{n-2}E_2^* ) > 3\log{n} - \alpha/8,
	\end{equation}
we have that 
	\eqn{
	\label{eqn:assump}
	(X_{\sss(1)},X_{\sss (2)})\stackrel{d}{=}\left(\frac{n}{n-2}E_1^*,\frac{n}{n-2}E_2^*\right).
	}
The probability of the event in \eqref{eqn:unless} is bounded by $2\e^{\alpha/8}/n^3$.	
Since $n^2 \e^{\alpha/8}{n^3}\to 0 $, we can ignore the contribution of this \ch{in the proof of} 
Proposition \ref{prop-small-min-edge-weight} and assume \eqref{eqn:assump}. 

\ch{Let} $V_1$ be the closest vertex to $1$, at distance $X_{\sss(1)}$ (respectively $V_2$ at distance $X_{\sss(2)}$ from vertex $2$). The rest of the \ch{smallest-weight} path has the same \ch{distribution} as the \ch{smallest-weight} path between 2 sets $A = \set{1, V_1}$ and $B=\set{2, V_2}$ in $\sS_n$.  \ch{Lemma \ref{lem-dist-sets-vert} thus implies that} 
	\eqn{
	d_w(i,j)=X_{\sss(1)}+X_{\sss (2)}+\sum_{k=2}^{N-1} \frac{nE_k^\prime}{k(n-k)},
	}
where $N=N_1\wedge N_2$ and $(N_1,N_2)$ is a uniform pair of distinct vertices from $[n]\setminus \set{1,2}$ and $(E_k^\prime)_{k\geq 1}$ are mean one exponential random variables.  Writing $S_N=\sum_{k=2}^{N-1} 
\frac{nE_k}{k(n-k)}$, we get 
	\eqn{
	\E[R_n(\alpha)]\leq n^2\pr\Big(S_{N}\geq 3\log{n}-X_{\sss(1)}-X_{\sss (2)}-\alpha/8,  X_{\sss(1)}\leq \log{n}-\alpha,
	X_{\sss(2)}\leq \log{n}+\alpha/2\Big).
	}
Thus,
	\eqn{
	\E[R_n(\alpha)]\leq n^2\int_0^{\log{n}-\alpha}\int_0^{\log{n}+\alpha/2}
	\e^{-(x+y)(n-2)/n}\pr\Big(S_{N}\geq 3\log{n}-x-y-\alpha/8\Big)dxdy.
	}
To complete the proof, we study the tail behavior of the random variable $S_N$. 

\begin{lem}[Tail behavior for random sums]
\label{lem-tail_SN}
For any constant $a<2$, there exists a $C=C_a$ such that for every $x\geq 0$,
	\eqn{
	\pr(S_N\geq \log{n}+x)\leq C\e^{-ax}.
	}
\end{lem}

\proof We compute the moment generating function of $S_N$ as
	\eqan{
	M_{S_N}(t)&=\sum_{j=2}^{n-2} \pr(N=j)\E[\e^{t S_j}]
	=\sum_{j=2}^{n-2} \pr(N=j) \prod_{k=2}^{j-1} \frac{k(n-k)}{k(n-k)-tn}\\
	&=\sum_{j=2}^{n-2} \pr(N=j)  \e^{-\sum_{k=2}^{j-1} \log(1-\frac{tn}{k(n-k)})}.\nn
	}
Thus,
	\eqan{
	\pr(S_N\geq \log{n}+x)&\leq \e^{-t(\log{n}+x)}M_{S_N}(t) \notag \\ 
	&\leq \e^{-t(\log{n}+x)} \sum_{j=2}^{n-2} \pr(N=j)  \e^{-\sum_{k=2}^{j-1} \log(1-\frac{tn}{k(n-k)})}.
	}
Take $t=a<2$ \ch{ and note that then $tn/[k(n-k)]<1$ since $k, n-k\geq 2$. Therefore, we can Taylor expand}
	\eqn{
	\log\left(1-\frac{tn}{k(n-k)}\right)\leq \frac{tn}{k(n-k)}+O(\frac{n^2}{[k(n-k)]^2}),
	}
Using that 
	\[
	\frac{n}{k(n-k)}=\frac{1}{k}+\frac{1}{n-k},
	\]
we arrive at
	\eqan{
	\pr(S_N\geq \log{n}+x)&\leq \e^{-t(\log{n}+x)}M_{S_N}(t)
	\leq C\e^{-a(\log{n}+x)} \sum_{j=2}^{n-2} \pr(N=j)  \e^{a\sum_{k=2}^{j-1} [\frac{1}{k}+\frac{1}{n-k}]}\nn\\
	&\leq C\e^{-ax} \sum_{j=2}^{n-2} \pr(N=j)  \e^{a[\log{(j/n)}-\log{(1-j/n)}]}\nn\\
	&=C\e^{-ax} \E\Big[\Big(\frac{N/n}{1-N/n}\Big)^a\Big].\nn
	}
Note that $\pr(N=j)=\frac{2(n-j)}{(n-2)(n-3)}$, so that, by dominated convergence,
	\eqn{
	\E\Big[\Big(\frac{N/n}{1-N/n}\Big)^a\Big]=\sum_{j=2}^{n-2} \frac{2(n-j)}{(n-2)(n-3)}
	\Big(\frac{j/n}{1-j/n}\Big)^a\to \int_0^1 \frac{u^a}{(1-u)^a} 2(1-u)du<\infty,
	}
whenever $a<2$.
\qed
%Write
%	\[
%	\frac{n}{k(n-k)}=\frac{1}{k}+\frac{1}{n-k},
%	\]
%and use that
%	\eqan{
%	\pr(\sum_{k=2}^{n-1} \frac{nE_k}{k(n-k)}\geq a)
%	&\leq 
%	\pr(\sum_{k=2}^{n-1} \frac{E_k}{k}\geq a/2)+\pr(\sum_{k=2}^{n-1} \frac{E_k}{n-k}\geq a/2)\\
%	&\leq 2\pr(\sum_{k=1}^{n-2} \frac{E_k}{k}\geq a/2).\nn
%	}
%Further,
%	\eqn{
%	\pr(\sum_{k=1}^{n-2} \frac{E_k}{k}\geq a)=\pr(\max_{i=1}^{n-1}E_i\geq a)
%	=1-(1-\e^{-a})^{n-1}\leq n\e^{-a}.
%	}
%Conditioning on $X_{\sss(1)}$ and applying this to $a=3\log{n}-X_{\sss(1)}-\alpha/2$ yields
%	\eqn{
%	\E[R_n(\alpha)]\leq 2n^2 
%	\E\Big[\indic{X_{\sss(1)}\leq \log{n}-\alpha} \e^{-(3\log{n}-X_{\sss(1)}-\alpha/2)/2)}\Big].
%	}
%	
%\qed
\bigskip

\noindent \ch{By Lemma \ref{lem-tail_SN}, with $a=3/2$,}
	\eqan{
	\E[R_n(\alpha)]&\leq Cn^2\int_0^{\log{n}-\alpha}\int_0^{\log{n}+\alpha/2}
	\e^{-(x+y)}\e^{-a(2\log{n}-x-y-\alpha/8)}dxdy\\
	&=Cn^{2a}\int_0^{\log{n}-\alpha}\int_0^{\log{n}+\alpha/2}
	\e^{(a-1)(x+y)}\e^{\alpha/8}dxdy\leq C\e^{-\alpha+\alpha/2+\alpha/4}\nn\\
	&= C\e^{-(a-1)\alpha/2+a\alpha/8}\leq C\e^{-\alpha/16}.\nn
	}
This completes the proof of Proposition \ref{prop-large-min-edge-weight}.
\qed
\vskip1cm
~

\noindent
\noindent
\subsection{{The limiting random variable}}
In this section, we prove the finiteness of the random variable $\Xi=\max_{s<t} (Y_s+Y_t-\Lambda_{st})$ in \eqref{eqn:xi-def} which Theorem \ref{thm:main} asserts is the limit of the re-centered diameter.
In the following lemma, we give an alternate expression for its distribution:

\begin{lem}[The limiting random variable]
\label{lem-lim-var}
Let $Q=\e^{-\Xi}$. Then, 
	\eqn{
	Q=\min_{s<t} \frac{S_sS_t}{E_{st}^\prime},
	} 
where $S_s=\sum_{i=1}^s E_i^\prime$ and
$(E_i^\prime)_{i\geq 1}$ and $(E_{st}^\prime)_{s<t}$ are i.i.d.\ exponential random variables with \ch{mean} 1.
In particular, for every $x>0$,
	\eqn{
	\label{tail-Q}
	\pr(Q>x)=\E\Big[\prod_{1\leq s<t} \big(1-\e^{-S_sS_t/x}\big)\Big],
	}
and $\pr(Q>x)\in (0,1)$ for every $x>0$.
\end{lem}

\proof We note that we can write $-\Lambda_{st}=\log(E_{st}^\prime)$ and $Y_s=-\log(S_s)$.
Indeed, the point process $(\e^{-Y_s})_{s\geq 1}$ is a standard Poisson process.
Thus, 	
	\eqn{
	\e^{-\Xi}\stackrel{d}{=}\min_{s<t} \e^{\log(S_s)+\log(S_t)-\log(E_{st}^\prime)}=Q.
	}
Equation \eqref{tail-Q} immediately follows. To prove that $\pr(Q>x)\in (0,1)$ for every $x>0$,
we note that $\pr(Q>x)<1$ follows immediately from \eqref{tail-Q} since each of the terms in the product
is $<1$ a.s. To show that $\pr(Q>x)>0,$ we first note that 
	\eqan{
	\pr(Q>x)&\ch{\geq}\E\Big[\prod_{1\leq s<t} \big(1-\e^{-S_sS_t/x}\big)\indic{S_1>1}\Big] \notag\\
	&=\E\Big[\prod_{1\leq s<t} \big(1-\e^{-S_sS_t/x}\big)\mid S_1>1\Big]\pr(S_1>1). \label{tail-Q-1}
	}
We compute that $\pr(S_1>1)=1/\e$, and observe that by the memoryless property of the exponential
random variable $S_1$, conditionally on $S_1>1$, the distribution of $(S_t)_{t\geq 1}$ is equal to 
$(S_t+1)_{t\geq 1}$. Thus, 
	\eqan{
	\label{tail-Q-2}
	\pr(Q>x)& \geq \e^{-1}\E\Big[\prod_{1\leq s<t} \big(1-\e^{-(S_s+1)(S_t+1)/x}\big)\Big]\nn\\
	&\geq \e^{-1}\exp\Big(\sum_{1\leq s<t} \E\Big[\log\big(1-\e^{-(S_s+1)(S_t+1)/x}\big)\Big].
	}
Next, we compute, using Fubini,
	\eqan{
	&\sum_{1\leq s<t} \E\Big[\log\big(1-\e^{-(S_s+1)(S_t+1)/x}\big)\Big]\\
	&\quad=\sum_{1\leq s<t} \int_0^{\infty} du \int_0^{\infty} dv 
	\frac{u^{s-1}}{(s-1)!}\frac{v^{t-s-1}}{(t-s-1)!}\e^{-(u+v)}\log\big(1-\e^{-(u+1)(v+1)/x}\big)\nn\\
	&\quad=\int_0^{\infty} du \int_0^{\infty} dv \sum_{1\leq s<t}
	\frac{u^{s-1}}{(s-1)!}\frac{v^{t-s-1}}{(t-s-1)!}\e^{-(u+v)}\log\big(1-\e^{-(u+1)(v+1)/x}\big)\nn\\
	&\quad=\int_0^{\infty}\int_0^{\infty}\log\big(1-\e^{-(u+1)(v+1)/x}\big)dudv<\infty.\nn
	}
This completes the proof.
\qed
\vskip1cm
~

\noindent
\subsection{{The limiting maximization problem}}
\label{sec:lim-max-prob}
In this section, we combine the various ingredients proved in the previous sections to prove the distributional convergence in Theorem \ref{thm:main}. We defer the proof of the convergence of moments to the next section.  By Proposition \ref{prop-large-min-edge-weight} and whp for large $\alpha$, 
$N_n(\alpha)\geq 2$. By Proposition \ref{prop-dist-large-min-edge-weight}, 	
	\eqn{
	\diam_w(K_n)-3\log{n}\geq d_w(V_1,V_2)-3\log{n}\weakc -2\alpha+\Lambda_1+\Lambda_2-\Lambda_{12}.
	}
As a result, $\diam_w(K_n)-3\log{n}\geq -K$ whp when $K>0$ is sufficiently large. Therefore,
also using Proposition \ref{prop-small-min-edge-weight}, whp for $\alpha$ sufficiently large,
	\eqn{
	\diam_w(K_n)=\max_{s<t\leq N_n(\alpha)}d_w(V_s,V_t).
	}
We note that, again using Proposition \ref{prop-dist-large-min-edge-weight} 
and Proposition \ref{prop-large-min-edge-weight},
	\eqn{
	\label{fin-alpha-conv}
	\max_{s<t\leq N_n(\alpha)} d_w(V_s,V_t)-3\log{n}
	\weakc \max_{s<t\leq N(\alpha)} (\Lambda_s+\Lambda_t-\Lambda_{st}-2\alpha),
	}
where $N(\alpha)$ is a Poisson random variable with mean $\e^{\alpha}$ and the Gumbel variables
are independent of $N(\alpha)$. As a result,
	\eqn{
	\diam_w(\ch{\cK_n})-3\log{n}\weakc \Xi^*,
	}
where $\Xi^*$  is the distributional limit as $\alpha\to \infty$ \ch{of the right-hand side of \eqref{fin-alpha-conv}, i.e.,}	
	\eqn{
	\max_{s<t\leq N(\alpha)} (\Lambda_s+\Lambda_t-\Lambda_{st}-2\alpha)\weakc \Xi.
	}
We show that this weak limit exists and that $\Xi^* = \Xi$ defined in \eqref{eqn:xi-def}. 

\begin{prop}[The limiting variable $\Xi$]
\label{prop-lim-var}
As $\alpha\to \infty$, 
	\eqn{
	\max_{s<t\leq N(\alpha)} (\Lambda_s+\Lambda_t-\Lambda_{st}-2\alpha)\weakc \Xi,
	}
where $\Xi$ is defined in \eqref{eqn:xi-def}.
\end{prop}
\proof 

As $\alpha\to \infty$,
	\eqn{
	\e^{-\alpha}N(\alpha)\Probc 1.
	}
Therefore, it suffices to prove that
	\eqn{
\Xi_\alpha:=\max_{s<t\leq \e^{\alpha}} (\Lambda_s+\Lambda_t-\Lambda_{st}-2\alpha)\weakc \Xi,
	}

Recall from Section \ref{sec:res}, the Poisson point process $\cP = \crh{(Y_s)_{s\geq 1}}$ with intensity measure given by the density function $\lambda(y) = \e^{-y}$. Also recall from \eqref{eqn:xi-def} that we defined $\Xi$ as 
	\[
	\crh{\Xi:= \max_{s< t}(Y_s+ Y_t - \Lambda_{st}).}
	\] 
For any fixed $A > 0$, let $\cP(A)$ denote $\cP$ restricted to the interval $[-A, \infty)$. Write 
	\[
	\crh{\Xi(A):= \max_{s< t\colon Y_s, Y_t \in \cP(A)}(Y_s+ Y_t - \Lambda_{st}).} 
	\]
Thus, $\crh{\Xi(A)}$ is the maximum of corresponding pairs $(s,t)$ whose point process values \crh{satisfy} $\crh{Y_s, Y_t \geq -A}$. Intuitively, \crh{one would expect that}  $\Xi = \Xi(A)$ for large $A$. \crh{We now make his intuition precise.} Define 
	\[
	\crh{\cR^{\sss(1)}(A):= \max_{s<t\colon Y_s, Y_t \leq  -A} (Y_s + Y_t - \Lambda_{st}),}
	\]
and, for $A < B$, let
	\[
	\crh{\cR^{\sss(2)}(A, B) := \max_{s<t\colon Y_s \geq  -A, Y_t \leq -(A+B)} (Y_s + Y_t - \Lambda_{st}).}
	\]
\crh{The random variable} $\cR^{\sss(1)}(A)$ is the supremum between pairs $(s,t)$ such that $Y_s , Y_t\leq -A$ while $\cR^{\sss(2)}(A,B)$ corresponds to supremum between pairs of points $(s,t)$ such that $Y_s > -A$ but $Y_t < -(A+B)$. Note that, for any $z$,
	\begin{equation}
	\label{eqn:set-xi-xia}
	\set{\Xi = \Xi(A+B)} \supseteq \set{\Xi(A) > z, \cR^{\sss(1)}(A) < z, \cR^{\sss(2)}(A, B) < z}.
	\end{equation} 

Consider the point process 
	\[
	\cP_\alpha^* = \sum_{s=1}^{\e^\alpha} \delta\set{\Lambda_s - \alpha}.
	\]
When arranged in increasing order, write this point process as $Y_1(\alpha) > Y_2(\alpha) > \cdots$. Standard extreme value theory implies that 
	\begin{equation}
	\label{eqn:cpalpha-cp}
	\cP_\alpha^* \weakc \cP \qquad \mbox{ as } \alpha \to \infty,
	\end{equation}
where $\weakc$ denotes convergence in distribution in the space of point measures on $\bR$ equipped with the vague topology.  Define, analogously to $\Xi(A), \cR^{\sss(1)}(A), \cR^{\sss(2)}(A, B) $, the random variables $\Xi_\alpha(A), \cR^{\sss(1)}_\alpha(A), \cR^{\sss(2)}_\alpha(A, B) $,  \crh{i.e.,}
	\[
	\crh{\Xi_\alpha(A):= \max_{s<t\colon Y_s(\alpha), Y_t(\alpha) \in \cP_A(\alpha)}(Y_s(\alpha) + Y_t(\alpha) - \Lambda_{st}).}
	\]
where $\cP_\alpha(A)$ is the point process $\cP_\alpha$ restricted to the interval $[-A, \infty)$. Similarly define $\cR_\alpha^{\sss(1)}(A),\cR^{\sss(2)}_\alpha(A)$. As before, for any $z$,
	\begin{equation}
	\label{eqn:set-xalpha-xalpha-A}
	\set{\Xi = \Xi(A+B)} \supseteq \set{\Xi(A) > z, \cR^{\sss(1)}(A) < z, \cR^{\sss(2)}(A, B) < z}
	\end{equation}

The weak convergence in \eqref{eqn:cpalpha-cp} immediately implies that, for any fixed $A$, 
	\begin{equation}
	\label{eqn:xialpha-xi-A}
	\Xi_\alpha(A) \weakc \Xi(A) \qquad \mbox{ as } \alpha \to\infty
	\end{equation}
The following lemma formalizes the notion that for large $A$, $\Xi = \Xi(A)$ whp and, similarly, when $\alpha$ is large $\Xi_\alpha(A) = \Xi_\alpha$ \crh{whp. This is achieved by showing} that for large $A$, each of the random variables $\cR^{\sss(1)}(A), \cR^{\sss(1)}_\alpha(A)$, and, for each fixed $A$, for sufficiently large $B$, $\cR^{\sss(2)}(A,B), \cR^{\sss(2)}_\alpha(A,B)$ take large negative values.  Using \eqref{eqn:xialpha-xi-A}, 
\eqref{eqn:set-xi-xia} and \eqref{eqn:set-xalpha-xalpha-A} completes the proof of Proposition
\crh{\ref{prop-lim-var}}. 

\begin{lem}
	\label{lem:a-a}
	\begin{enumeratea}
		\item Fix $x\in \bR$. Then,
			\[
			\limsup_{A\to \infty} \pr(\cR^{\sss(1)}(A) > x) = 0.
			\] 
		Further, for each fixed $A$,
			\[
			\limsup_{B\to\infty} \pr(\cR^{\sss(2)}(A,B)> x) = 0.
			\]
		\item Fix $x\in \bR$. Then,
		 	\[
			\limsup_{A\to \infty}\limsup_{\alpha \to\infty} \pr(\cR^{\sss(1)}_{\crh{\alpha}}(A) > x) = 0.
			\] 
		Further, for each fixed $A$,
			\[
			\limsup_{B\to\infty} \limsup_{\alpha \to\infty} \pr(\cR^{\sss(2)}_{\crh{\alpha}}(A,B)> x) = 0.
			\]
	\end{enumeratea}
\end{lem}

\proof We \crh{start by proving part (a).} We start with $\cR^{\sss(1)}(A)$. To simplify notation, we also \crh{restrict ourselves to}  the case $x=0$. The general $x$ case is identical. 

Write 
	\[
	\cN^{\sss(1)}(A):= \# \set{(s,t): Y_s, Y_t < -A,  Y_s + Y_t- \Lambda_{st} \geq 0 }.
	\]
It is enough to show $\limsup_{A\to\infty}\E(\cN^{\sss(1)}(A)) = 0$. Conditioning on the point process $\cP$, we get 
	\[
	\E(\cN^{\sss(1)}(A)|\cP) = \sum_{(s,t), s< t, Y_s, Y_t < -A} \e^{-\e^{-(Y_s +Y_t)}}.
	\]
Fix $a> 1$. \crh{We use the fact that we can choose $A$ so large such that $\e^{-\e^{C+D}} < \e^{-aC} \e^{-aD}$ for all $C,D > A$. This leads to}
	\[
	\E(\cN^{\sss(1)}(A)|\cP) \leq \sum_{(s,t), s< t, Y_s, Y_t < -A} \e^{aY_s} \e^{aY_t}.
	\]
Since $\set{Y_s\in \cP: Y_s\leq -A}$ is just a Poisson point process on the interval $(-\infty, -A]$ with density $e^{-x}$, properties of Poisson processes \cite[Eqn 3.14]{kingman-book} implies that, as $A\to\infty$,
	\begin{align*}
	\E\Big(\sum_{\substack{(s,t), s< t,\\ Y_s, Y_t < -A}} \e^{aY_s} \e^{aY_t}\Big) 
	&=\frac{1}{2}\left(\int_{-\infty}^{-A} \e^{ax}\e^{-x} dx\right)^2\\ &= \frac{1}{2}\e^{-2(a-1)A}\to 0.
	\end{align*}
%\todo{Real equality? What about $s<t$?}
This shows that $\limsup_{A\to\infty}\E(\cN^{\sss(1)}(A)) = 0$ and thus completes the proof. 

Next fix $A$ and let us deal with $\cR^{\sss(2)}(A,B)$. Here we use the fact that $\cP(A)$ and $\cP^c(A+B) := \cP \setminus \cP^c(A+B)$ are independent Poisson point processes on the sets $[-A, \infty)$ and $(-\infty, -(A+B))$ with intensity measure with density $\lambda(y) = \e^{-y}$. We work conditional on $\cP(A)$. Fix a point $Y_s$ in $\cP(A)$. Then,
	\[
	\pr(\sup_{Y_t < -(A+B)} (Y_s + Y_t -\Lambda_{st}) < z |\cP(A)) 
	= \E\Big(\prod_{t\colon Y_t < -(A+B)} \Big(1-\e^{-\e^{-(Y_t -(z-Y_s))}}\Big) \Big).
	\]
The following lemma completes the proof:
\begin{lem}
\label{lem:sup-one}
	Fix any $z^*$ and $A$. Then 
	\[
	\lim_{B\to \infty}\E\Big(\prod_{t: Y_t < -(A+B)} \left(1-\e^{-\e^{-(Y_t -z^*)}}\right) \Big) \to 1.
	\]	
\end{lem}
\proof By the dominated convergence theorem, it is enough to show that, as $B\to\infty$,
	\[
	\prod_{t: Y_t < -(A+B)} \left(1-\e^{-\e^{-(Y_t -z^*)}}\right) \Probc 1.
	\]
Taking logarithms, this is equivalent to showing that, as $B\to\infty$,
	\[
	\sum_{t: Y_t < -(A+B)}\log{\left(1-\e^{-\e^{-(Y_t -z^*)}}\right)} \Probc 0.
	\]
In turn, this is equivalent to showing that, as $B\to\infty$,
	\[
	\sum_{t: Y_t < -(A+B)} \e^{-\e^{-(Y_t -z^*)}} \Probc 0.
	\]
By Campbell's theorem \cite{kingman-book}, 
	\begin{align*}
	\E(\sum_{t: Y_t < -(A+B)} \e^{-\e^{-(Y_t -z^*)}}) &= \int_{-\infty}^{-(A+B)} \e^{-\e^{-(y-z^*)}} \e^{-y} dy\\
	&= \e^{z^*} \e^{-\e^{ A+ B+z^*}} \to 0,
	\end{align*}
as $B\to\infty$. This completes the proof \crh{of part (a).} 

\crh{For part (b), we follow the proof of part (a). We highlight some of the differences only. We again start with $\cR^{\sss(1)}_{\alpha}(A)$ and again \crh{restrict ourselves to}  the case $x=0$. The general $x$ case is identical. 

Write 
	\[
	\cN^{\sss(1)}_{\alpha}(A):= \# \set{(s,t): Y_s(\alpha), Y_t(\alpha) < -A,  Y_s(\alpha) + Y_t(\alpha)- \Lambda_{st} \geq 0 }.
	\]
It is enough to show 
$\limsup_{A\to\infty}\limsup_{\alpha\rightarrow \infty}\E(\cN^{\sss(1)}_{\alpha}(A)) = 0$. Conditioning on the point process $\cP^*_{\alpha}$, we now get 
	\eqan{
	\E(\cN^{\sss(1)}_{\alpha}(A)|\cP^*_{\alpha}) 
	&= \sum_{(s,t), s< t, Y_s(\alpha), Y_t(\alpha) < -A} \e^{-\e^{-(Y_s(\alpha) +Y_t(\alpha))}}\\
	&=\sum_{1\leq s<t\leq \e^{\alpha}} \indic{\Lambda_s, \Lambda_t< -A+\alpha} \e^{-\e^{-(\Lambda_s-\alpha)-
	(\Lambda_t-\alpha)}}.\nn
	}
Now taking expectations and using that $\Lambda_s,\Lambda_t$ are independent for $s<t$ leads to
	\eqn{
	\E(\cN^{\sss(1)}_{\alpha}(A))\leq 
	\int_{-\infty}^{-A+\alpha}\int_{-\infty}^{-A+\alpha} \e^{-(u-\alpha)}\e^{-\e^{-u}}
	\e^{-(v-\alpha)}\e^{-\e^{-v}}\e^{-\e^{-(u-\alpha)-
	(v-\alpha)}}dudv.
	}
This integral can be bounded by
	\eqn{
	\E(\cN^{\sss(1)}_{\alpha}(A))\leq 
	\int_{-\infty}^{-A}\int_{-\infty}^{-A} \e^{-u}
	\e^{-v}\e^{-\e^{-u-
	v}}dudv,
	}
which is independent of $\alpha$ and converges to 0 as $A\rightarrow \infty$. The proof for 
$\cR^{\sss(2)}_{\crh{\alpha}}(A,B)$ is similar and will be omitted.
}
\qed

\vskip1cm
~

\subsection{{Convergence of moments}}
Recall that \crh{$C_{ij} = d_w(i,j)$.} \crh{We need to show 
	\[
	\E[\max_{i,j\in [n]} C_{ij}] - 3\log{n} \to \E[\Xi],
	\qquad
	{\rm Var}(\max_{i,j\in [n]} C_{ij})\to {\rm Var}(\Xi).\]
	}
%\todo{Notation clash, $C_{ij}$ is now a factor $n$ larger than in Theorem \ref{thm:main}.}

\csb{Since we have already shown convergence in distribution,} by uniform integrability for any $p\geq 1$, to prove that
	\eqn{
	\E\Big[\big(\max_{i,j\in [n]} C_{ij} - 3\log{n}\big)^p\Big]\to \E[\Xi^p],
	}
it suffices to prove that, \ch{for some integer $q$ with $q>p/2$,}
	\eqn{
	\E\Big[\big(\max_{i,j\in [n]} C_{ij} - 3\log{n}\big)^{2q}\Big]=O(1).
	}

\ch{Combined with convergence in distribution, this} implies convergence of the moments as well 
as existence of the moments of the limit random variable $\Xi$. Note that
	\eqn{
	\label{split-moments}
	\E\Big[\big(\max_{i,j\in [n]} C_{ij} - 3\log{n}\big)^{2q}\Big]=
	\E\Big[\big(\max_{i,j\in [n]} C_{ij} - 3\log{n}\big)^{2q}_+\Big]+
	\E\Big[\big(\max_{i,j\in [n]} C_{ij} - 3\log{n}\big)^{2q}_-\Big].
	}
\ch{We start by analyzing the first term on the right-hand side of \eqref{split-moments} by deriving an upper bound
on $\max_{i,j\in [n]} C_{ij} - 3\log{n}$, and then prove a lower bound on  $\max_{i,j\in [n]} C_{ij} - 3\log{n}$ 
to obtain a bound on the second term on the right-hand side of \eqref{split-moments}.}
\bigskip

\noindent{\bf Upper bound:}	
Let us analyze the first term and show that 
	\[
	\E\Big[\big(\max_{i,j\in [n]} C_{ij} - 3\log{n}\big)^{2q}_+\Big]= O(1).
	\]
To prove this assertion, it is enough to show that there exist $N,\alpha$ such that for all large 
$n > N$ and $x\geq \alpha$, \ch{the random variable $\max_{i,j\in [n]} C_{ij} - 3\log{n}$ 
has exponential upper tails} in the sense that there exist constants $\kappa_1, \kappa_2> 0$ 
(independent of $x$) such that 
	\begin{equation}
	\label{eqn:exponent-tail-bd}
	\pr(\max_{i,j\in [n]} C_{ij} - 3\log{n} > x) \leq \kappa_1\e^{-\kappa_2 x}.
	\end{equation}
Now note that 
	\eqn{
	\label{split-event-upper-tails}
	\ind\set{\max_{i,j\in [n]} C_{ij} - 3\log{n} > x} 
	\leq \ind\set{\max_{i\in [n]} X_{\sss(i)} > \log{n}+ 4x}
	+ R_n^{\sss(1)}(x) 
	+ R_n^{\sss(2)}(x).
	}
Here $R_n^{\sss(1)}(x) =R_n(8x)$ as in \eqref{eqn:rn-def}, i.e.,
	\[
	R_n^{\sss(1)}(x)
	=\#\set{i,j\in [n]\colon X_{\sss(i)}\leq \log{n}-8x, X_{\sss(j)}\leq \log{n}+4x, d_w(i,j)\geq 3\log{n}-x},
	\] 
while 
	\[
	R_n^{\sss(2)}(x)
	:= \#\set{(i,j)\colon X_{\sss(i)}> \log{n}-8x, X_{\sss(j)}> \log{n}-8x, d_w(i,j) > 3\log{n}+x}.
	\]
Recall that for any $\alpha\in \bR$,  $N_n(\alpha)$ denotes the number of vertices $i$ with $X_{\sss(i)} \geq \log{n} -\alpha$.  \ch{For the first term in \eqref{split-event-upper-tails}, since $\pr(\max_{i\in [n]} X_{\sss(i)} > \log{n}+ 4x) = \pr(N_n(-4x) \geq 1)$, the Poisson approximation in Proposition \ref{prop-large-min-edge-weight} implies that}
	\begin{align}
	\pr(\max_{i\in [n]} X_{\sss(i)} > \log{n}+ 4x) &	
	\leq \frac{2(1+o(1))\e^{-4x} \log{n}}{n}+ (1-\e^{-\e^{-4x}})\nn \\ 
	&\leq (1+o(1))\e^{-4x}.	
	\label{eqn:ind-1} 
	\end{align}
Further, by Proposition \ref{prop-small-min-edge-weight} for $n$ large enough
	\begin{equation}
	\label{eqn:rn1-bd}
	\E(R_n^{\sss(1)}(x)) \leq C\e^{-x/2}. 
	\end{equation}
We are left \ch{to analyze} $R_n^{\sss(2)}(x)$. Arguing as in the proof of Proposition \ref{prop-small-min-edge-weight},
	\[
	\E(R_n^{\sss(2)}(x)) \leq \E(N_n^2(-8x)) \pr(nd_w(1,2)> \log{n}+17x),
	\] 
where $d_w(1,2)$ is the distance between vertices $1,2$ in $\sS_n = \set{\cK_n, (E_e)_{e\in \cE_n}}$. Since 
	\[
	d_w(1,2)\stackrel{d}{=} \sum_{k=1}^{N} \frac{E_j}{k(n-k)},
	\]
where $N$ is uniform on $[n-1]$ independent of $(E_j)_{j\in[n-1]}$ which are mean $n$ exponential 
random variables. Thus, by Markov's inequality,  for any $\alpha> 0$ 
 	\[\pr(d_w(1,2) - \log{n} > 17x) 
	\leq \e^{-17\alpha x}\sum_{j=1}^{n-1}\frac{1}{n-1} 
	\exp\biggl(\alpha \biggl[\log\frac{j}{n} - \log{\biggl(1-\frac{j}{n}\biggr)}\biggr]\biggr).
	\]
Letting $\beta = 1-\eps$ with $\eps> 0$ small but independent of $x,n$, we finally get 	
	\[
	\pr(d_w(1,2) - \log{n} > 17x) \leq (1+o(1)) \e^{-17\alpha x} \E\left( \left[\frac{U}{1-U}\right]^{1-\eps}  \right),
	\]
where $U\sim U[0,1]$. We need to now bound $\E(N_n^2(-8x))$. Write $N_n(-8x) = \sum_{i=1}^n Z_i$ where $Z_i = \ind\set{X_{\sss(i)}\geq \log{n}+8x}$. By Proposition \ref{prop-large-min-edge-weight}, 
$\E(N_n(-8x))\leq 2\e^{8x}$. Further, 
	\[
	\var(N_n(-8x)) \leq 2\e^{8x} + n(n-1)\pr(Z_1=1)[\pr(Z_2=1|Z_1=1) - \pr(Z_2=1)].
	\] 
Given $Z_1 =1$, the edge weights $\ch{(E_{2,i})_{i\neq 2}}$ have \ch{the same distribution as} $\big(\log{n}-8x+ E_{2,1}, (E_{2,j})_{j\neq 1,2}\big)$. Thus,
	\[
	\pr(Z_2 = 1|Z_1= 1) = \pr(\min_{j\geq 2} E_{2,j} > \log{n} - 8x) = \exp\left(-\frac{n-2}{n}(\log{n} -8x)\right) . 
	\]
Combining this, we get that $\var(N_n(-8x)) \leq 4 \e^{8x}$ so that $\E([N_n(-8x)]^2) \leq 16 \e^{16x}$. This results in 
	\begin{equation}
	\label{eqn:rn2-bd}
	\E(R_n^{\sss(2)}\crh{(x)}) \leq (1+o(1)) 16 \E\left( \left[\frac{U}{1-U}\right]^{1-\eps}  \right) \e^{-(1-17\eps)x}.
	\end{equation}
Combining \eqref{eqn:ind-1}, \eqref{eqn:rn1-bd} and \eqref{eqn:rn2-bd} completes the proof of the asserted exponential tail bound in \eqref{eqn:exponent-tail-bd} and completes the proof of the upper bound. 
\ \\ \ \\ 

\noindent {\bf Lower bound:} Let us now show that
	\[
	\E\Big[\big(\max_{i,j\in [n]} C_{ij} - 3\log{n}\big)^{2q}_-\Big] = O(1).
	\]
Recall that $V_1, V_2$ denote the vertices with the largest and second largest $X_{\sss(i)}$ values. Further 
	\[
	\max_{i,j\in [n]} C_{ij} - 3\log{n} \geq_{st} \crh{(X_{\sss (V_1)} - \log{n})_{-}
	+ (X_{\sss (V_2)} - \log{n})_{-} + (nd_w(1,2)-\log{n}),}
	\]
where $d_w(1,2)$ is independent of $X_{\sss (V_i)}$ with the same distribution as the length of the optimal path between $1,2$ in $\sS_n$ and $\geq_{st}$ denotes stochastic domination. By H\"older's inequality
	\eqan{
	\E\Big[\big(\max_{i,j\in [n]} C_{ij} - 3\log{n}\big)^{2q}_-\Big] 
	&\leq 3^{2q} \Big(\E\left(\left[X_{\sss (V_1)}-\log{n}\right]_{-}^{2q}\right) 
	+ \E\left(\left[X_{\sss (V_2)}-\log{n}\right]_{-}^{2q}\right)\nn\\
	&\qquad + \E\left(\left[d_w(1,2)-\log{n}\right]^{2q}\right)\Big).
	}
By \cite[Proof of Theorem 3.3]{janson-cm-gp}
	\[
	\E\left(\left[d_w(1,2)-\log{n}\right]^{2q}\right) = O(1).
	\]
Further, $\E\left(\left[X_{\sss (V_1)}-\log{n}\right]_{-}^{2q}\right) \leq \E\left(\left[X_{\sss (V_2)}-\log{n}\right]_{-}^{2q}\right)$. 
Using the identity 
	\[
	\E(Y^{2q}) = (2q-1) \int_0^\infty y^{2q-1} \pr(Y> y) dy,
	\] 
for any \ch{non-negative} random variable $Y$ and $\left[X_{\sss (V_2)}-\log{n}\right]_{-}^{2q} \leq (\log{n})^{2q}$, it is enough to show for some $0<\eps< 1$ small enough
	\begin{equation}
	\label{eqn:bd-lwb-enough}
	\pr(\log{n}- X_{\sss (V_2)} \geq x )\leq \left\{\begin{array}{ll}
	2{\e^{-(1-\eps)\e^{x}}} + 2\frac{\e^{2x} \log{n}}{n}, \qquad & x< (1-\eps)\log{n}/2, \\
	\e^{-n^{1/3}} + \frac{\log{n}}{\sqrt{n^{1/3}}} \qquad & x\in [(1-\eps)\log{n}/2,\log{n}]. 
	\end{array}\right.
	\end{equation}
The first line follows from the Poisson approximation result Proposition \ref{prop-large-min-edge-weight} since $\pr(\log{n}- X_{\sss (V_2)} \geq x ) = \pr(N_n(x) \leq 1)$. To prove the second line consider the case where $x=(1-\eps)\log{n}/2$. Fix a set $\cA\subseteq [n]$ with size $|\cA| = n^{1/3}$. For each vertex $v\in \cA$, define 
	\[
	X^*_{\sss(v:[n]\setminus \cA) } = \min_{j\in [n]\setminus\cA} E_{v,j}.
	\] 
Then $\ch{(X^*_{\sss(v:[n]\setminus \cA)})_{v\in \cA}}$ is a collection of $n^{1/3}$ independent exponential mean $n/(n-n^{1/3})$ random variables. Define $N_n^* = \sum_{v\in \cA} \ind\set{X^*_{\sss(v:[n]\setminus \cA)} \geq (1+\eps)\log{n}/2}$. Then one can check that 
	\begin{equation}
	\label{eqn:v2-bd}
	\ind\set{X_{\sss (V_2)} <  (1+\eps) \log{n}/2 } 
	\leq \ind\set{\min_{i,j\in \cA} E_{i,j}< (1+\eps)\log{n}/2 } + \ind\set{N_n^* \leq 1},	
	\end{equation}
\ch{since $\min_{i,j\in \cA} E_{i,j}< (1+\eps)\log{n}/2$ and $X_{\sss (V_2)} <  (1+\eps) \log{n}/2$ implies that $N_n^* \leq 1$.}
Now note that 
	\[
	N_n^* \sim \mathrm{Bin}\left(n^{1/3},1-\exp\Big(-\frac{(n-n^{1/3})}{n} (1+\eps)\log{n}/2\Big) \right),
	\]
while $\min_{i,j\in \cA} E_{i,j}$ has \ch{an exponential distribution with} rate $\ch{n^{1/3}(n^{1/3}-1)/(2n)}$ since the number of edges in $\cA$ is $\ch{n^{1/3}(n^{1/3}-1)/2}$. Further,
	\[
	\ch{1-\exp\left(-\frac{(n-n^{1/3})}{n}\log{n}/2\right) \geq 1-\frac{1}{n^{1/3}}.}
	\]
Taking expectations in \eqref{eqn:v2-bd} completes the proof of \eqref{eqn:bd-lwb-enough} and thus the proof of the lower bound. This completes the proof of the main result.  \qed 

%To show that $\pr(Q>x)>0,$ we note that by the FKG inequality and the fact that 
%	\[
%	1-\e^{-S_sS_t/x}
%	\]
%is a montone functional of $(E_s)_{s\geq 1}$, we obtain
%	\eqn{
%	\pr(Q>x)=\E\Big[\prod_{1\leq s<t<\infty} \big(1-\e^{-S_sS_t/x}\big)\Big]
%	\geq \prod_{1\leq s<t<\infty} \E[1-\e^{-S_sS_t/x}].
%	}
%Now, 
%	\eqn{
%	\E[1-\e^{-S_sS_t/x}]=1-\E[\e^{-S_s^2/x}\e^{-S_s(S_t-S_s)/x}]
%	=1-\E\big[\e^{-S_s^2/x}\big(\frac{1}{1+S_s/x}\big)^{t-s}\big].
%	}
%Further,
%	\eqan{
%	\E\big[\e^{-S_s^2/x}\big(\frac{1}{1+S_s/x}\big)^{t-s}\big]
%	&=\int_0^{\infty} \frac{u^{s-1}}{(s-1)!}\e^{-u} \e^{-u^2/x}\big(\frac{1}{1+u/x}\big)^{t-s}du\\
%	&\leq \int_0^{\infty} \frac{u^{s-1}}{(s-1)!}\e^{-u} \e^{-u^2/x}\big(\frac{1}{1+u/x}\big)^{t-s}du.\nn
%	}
%For $s=1$,
%	\eqan{
%	\E\big[\e^{-S_1^2/x}\big(\frac{1}{1+S_1/x}\big)^{t-1}\big]
%	&\leq \int_0^{\infty} \e^{-u}\big(\frac{1}{1+u/x}\big)^{t-1}du\\
%	&=\int_0^{\infty} \e^{-u-(t-1)\log(1+u/x)}du\sim x/(t-1),\nn
%	}
%which is unfortunately not summable over $t$... Mmmm. Is
%	\eqn{
%	\E\Big[\prod_{1<t<\infty} \big(1-\e^{-S_1S_t/x}\big)\Big]>0?
%	}
%	
%		
%	
%\todo{Add simplest possible argument!}
%
%\qed
\vskip.3in
\noindent{\bf Acknowledgments.}  
The work of RvdH was supported in part by the Netherlands Organisation for Scientific Research (NWO).  \ch{The work of SB has been supported in part by NSF-DMS grant 1105581 and in part by an NWO Star grant. SB thanks the hospitality of {\sc Eurandom} where this work commenced in November 2012. We thank Julia Komj\'athy for a careful reading of an early version of the paper. }

%\listoftodos

% \bibliographystyle{plain}
% \bibliography{complete-graph}
% 
\def\cprime{$'$}

\end{document}